\documentclass[twoside]{aiml}

\usepackage{aimlmacro}

\usepackage{amsmath}
\usepackage{pgffor}
\usepackage{enumitem}
\usepackage{ifthen}
\usepackage{stmaryrd}
\usepackage{hyperref} \hypersetup{ colorlinks=true, linkcolor=blue, filecolor=blue,
citecolor = black,urlcolor=cyan, } \usepackage{diagbox} \usepackage[all]{xy}
\usepackage{cleveref}

\newcommand{\ttt}[1]{\texttt{#1}}

\renewcommand{\le}{\leqslant}
\renewcommand{\ge}{\geqslant}
\renewcommand{\models}{\vDash}

\newcommand{\md}{\mathrm{md}}

\newcommand{\NN}{\mathbb{N}}

\newcommand{\wV}[1][]{ \ifthenelse{\equal{#1}{}} {\ensuremath{\widehat{V}}}
{\ensuremath{\widehat{V_{#1}}}} } 
\newcommand{\tset}[1]{\llbracket #1 \rrbracket}

\newcommand{\bff}{\mathbf{F}} \newcommand{\bfg}{\mathbf{G}} \newcommand{\bbf}{\mathbb{F}}
\newcommand{\fv}{\mathsf{Fv}} \newcommand{\Prop}{\mathsf{Prop}}
\newcommand{\sfp}{\mathsf{p}} \newcommand{\KFr}{\mathsf{KFr}}
\newcommand{\ttq}{\mathtt{Q}} \newcommand{\ttat}{\mathtt{At}}
\newcommand{\ttbc}{\mathtt{Bc}} \newcommand{\ttr}{\mathtt{R}}
\newcommand{\sfa}{\mathsf{A}} \newcommand{\sfe}{\mathsf{E}}
\newcommand{\pqml}{\textsc{pqml}}
\newcommand{\npqml}{\textsc{npqml}}

\def\lastname{Ding, Li}

\begin{document}

\begin{frontmatter}
  \title{Some General Completeness Results for Propositionally Quantified Modal Logics}
  \author{Yifeng Ding, Yipu Li}
  \address{Department of Philosophy and Religious Studies, Peking University}
  
  \begin{abstract} 
We study the completeness problem for propositionally quantified modal logics on quantifiable general frames, where the admissible sets are the propositions the quantifiers can range over and expressible sets of worlds are admissible, and Kripke frames, where the quantifiers range over all sets of worlds. We show that any normal propositionally quantified modal logic containing all instances of the Barcan scheme is strongly complete with respect to the class of quantifiable general frames validating it. We also provide a sufficient condition for the truth of all formulas, possibly with quantifiers, to be preserved under passing from a quantifiable general frame to its underlying Kripke frame. This is reminiscent of both the idea of elementary submodel in model theory and the persistence concepts in propositional modal logic. The key to this condition is the concept of finite diversity (Fritz 2023), and with it, we show that if $\Theta$ is a set of Sahlqvist formulas whose class of Kripke frames has finite diversity, then the smallest normal propositionally quantified modal logic containing $\Theta$, Barcan, a formula stating the existence of world propositions, and a formula stating the definability of successor sets, is Kripke complete. As a special case, we have a simple finite axiomatization of the logic of Euclidean Kripke frames.
  \end{abstract}

\end{frontmatter}

\section{Introduction}

Propositionally quantified modal logics (\pqml s henceforth) are modal logics augmented with propositional quantifiers, a special kind of quantifiers that can be intuitively understood as capturing the quantification implicit in English sentences such as ``Everything Jane believes is false'' and ``It's likely that there's something I will never know''. One can also understand the expression ``For all the police knows, John is dead already'' as ``That John is dead already is compatible with everything the police knows'' and see that propositional quantification is involved.

Normal propositional modal logic has been studied fruitfully in relation to the possible
world-based Kripke frames. For propositional quantifiers, it is natural to consider them
as quantifying over sets of possible worlds, since in Kripke frames, propositional
variables are interpreted as sets of worlds, and propositional quantifiers bind these
variables.\footnote{This is obviously not the only way: another well known approach to
propositionally quantified formulas is to view them as uniform interpolants
\cite{pitts1992interpretation,visser1996bisimulations,Dagostino2000logical,bilkova2007uniform},
and semantically we also have the related \emph{bisimulation quantifiers}
\cite{ghilardi1995undefinability,French2006thesis,french2007idempotent,Dagostino2005axiomatization,steinsvold2020some,Dekker2024KD45}.
A detailed comparison is well beyond the scope of this paper, and we
only note that understanding propositional quantifiers as bisimulation quantifiers would
immediately lead to the \emph{atomless principle} $\forall p (\Diamond p \to \exists q
(\Diamond (p \land q) \land \Diamond (p \land \lnot q)))$, assuming a reasonably rich
model class. As we will see, one of the hallmarks of Kripke frames is that they validate
the opposite \emph{atomicity principle}, and quantifiable frames in general are silent on
this issue.} This immediately gives \pqml s a second-order flavor, and indeed they are
also known as second-order propositional modal logics when interpreted on Kripke frames.
We may also impose a distinction between sets of worlds that count as propositions
(possible value of propositional variables) and sets of worlds that do not count, thereby
adding a \emph{propositional domain} (also called the set of admissible sets) to Kripke
frames and obtaining general frames. General frames in which every formula under every
variable assignment expresses a proposition are called \emph{quantifiable} frames, and
they have the desirable logical property of validating instantiation reasonings such as
from ``Everything I believe is true'' and ``I believe that the Moon is made of cheese'' to
``The Moon is made of cheese''.

The theory of \pqml s based on Kripke frames and quantifiable frames has been studied
since the early days of modal logic (see \cite{Fritz2024} for a nice survey). The early
landmark paper is Fine's \cite{fine1970propositional} in which, among many other things,
it is shown that the \pqml\ consisting of formulas valid on the class of reflexive and
transitive Kripke frames is not recursively axiomatizable while the \pqml\ of the class of
Kripke frames with a universal relation is decidable and can be axiomatized by simply
adding to the modal logic $\mathtt{S5}$ the standard quantificational logic and an
atomicity principle stating that there is always a world proposition that is itself true
and entails every truth. Related questions of decidability, axiomatizability, and
expressivity under various kinds of frames were under sustained investigation
\cite{kremer1997complexity,kaminski1996expressive,antonelli2002representability,kuhn2004simple,ten2006expressivity,kuusisto2015second,belardinelli2018secondAIJ,kremer2018completeness,French2002AiML,zach2004decidability,Ding2021,bednarczyk2022whydoes}.
A recent breakthrough is by Fritz in \cite{fritz2022} where he established many new
results on the decidability and axiomatizability of \pqml s in a very general fashion.
Prior to this work, it was not even known whether the \pqml\ of Euclidean Kripke frames is
decidable, and the decidability of the \pqml\ of the Kripke frames validating
$\mathtt{KD45}$ was only established in \cite{Ding2021}.

In this paper, we focus on the question of when an axiomatically defined \pqml\ is
complete w.r.t. the quantifiable/Kripke frames it defines. In particular, we aim for a
counterpart of the celebrated Sahlqvist Completeness Theorem for \pqml s. For example, the
$\ttt{5}$ axiom $\Diamond p \to \Box\Diamond p$ defines the class of Euclidean Kripke
frames, and we know by Sahlqvist Completeness Theorem that $\ttt{5}$ is also sufficient to
axiomatize the normal modal logic of Euclidean Kripke frames. But is $\ttt{5}$ sufficient
also for the \pqml\ of Euclidean Kripke frames, even with the help of some atomicity
principles? A full generalization of the Sahlqvist Completeness Theorem was already
declared impossible by Fine's results, as the \pqml\ of the class of all Kripke frames is
not recursively axiomatizable but the class is trivially definable. Thus, we must build on
some general condition for axiomatizability, and \cite{fritz2022} provided an elegant one:
the finiteness of the \emph{diversity} of Kripke frames. Given a Kripke frame $\bbf$, two
worlds in it are called duplicates if the permutation switching them is an automorphism of
$\bbf$. Being duplicates is an equivalence relation, and the diversity of $\bbf$ is the
cardinality of the set of equivalence classes of this relation, while the \emph{diversity}
of a class $\mathsf{C}$ of Kripke frames is the supremum of the diversity of all
\emph{point-generated} subframes of frames in $\mathsf{C}$. Fritz \cite{fritz2022} shows
that if a class $\mathsf{C}$ of Kripke frames is defined by a finite set $\Theta$ of
formulas and has finite diversity, then the \pqml\ of $\mathsf{C}$ is decidable.
Decidability itself does not provide much information on whether there can be a simple and
intuitive axiomatization or if a given axiomatically defined logic is complete, but as we
mentioned above, Fine gave a simple axiomatization of the \pqml\ of Kripke frames with a
universal relation (equivalently, Kripke frames validating \ttt{S5}), and more recently,
for Kripke frames validating $\mathtt{KD45}$, it is shown \cite{Ding2021} that we only
need to add the quantificational axioms and rules, the Barcan scheme $\ttt{Bc}$, and the
atomicity principle for completeness. Can this result be generalized? Our finding is that
if $\Theta$ is a set of Sahlqvist formulas and the class $\mathsf{KFr}(\Theta)$ of Kripke
frames it defines has diversity $n$, then the normal \pqml\ axiomatized by $\Theta$, the
Barcan scheme $\ttt{Bc}$, an atomicity principle $\mathtt{At}^n$, and an axiom
$\mathtt{R}^n$ stating that the successor sets of propositions are propositions, is sound
and complete w.r.t.  $\mathsf{KFr}(\Theta)$. The last two axioms are parametrized by $n$
as they use iterated modalities to simulate the reflexive and transitive closure of the
primitive modality.  While we can replace the condition of $\Theta$ being Sahlqvist by
suitable technical conditions, we cannot drop it completely; and while for $\mathtt{KD45}$
$\mathtt{R}^n$ is redundant, for $\mathtt{K5}$ it is not.

Our method is based on saturated (witnessed) maximally consistent sets and the saturated
canonical general frame built from them. Fine \cite{fine1970propositional} claims that
this method can be used to show the completeness of $\mathtt{S5}$ with atomicity principle
w.r.t. Kripke frames with a universal relation but provided only an extremely terse
sketch. We will first show that if $\Lambda$ is a normal \pqml\ with the Barcan scheme
\ttt{Bc}, then its saturated canonical general frame is quantifiable and indeed validates
$\Lambda$. A corollary is that any normal $\pqml$ containing $\ttt{Bc}$ is sound and
complete w.r.t. the class of quantifiable frames it defines. Note that since quantifiable
frames are essentially multi-sorted first-order structures, \pqml s on them are no longer
second-order and we are not bound by the non-axiomatizability issues. This basic general
completeness result may have been realized by multiple scholars before, but a formal proof
for the fully general statement appears to be missing. 

For the main theorem, we start with the saturated canonical general frame and gradually
turn it into a Kripke frame that validates the original logic and satisfies a given
consistent set of formulas. The strategy is as follows: first (1) take a point-generated
general frame, then (2) keep only the worlds that are named by world propositions and
obtain what we call the \emph{atomic general frame}, and (3) finally show that we can
expand the propositional domain of the atomic general frame to the full powerset without
affecting the semantic value of any formula.  Step (2) also happens in some completeness
proofs of hybrid logics (see e.g. \cite{ten2004model}, Definition 5.2.6). The resulting
atomic general frame has the special property of being discrete and tense, and Sahlqvist
formulas are shown to be persistent over these general frames \cite{venema1993derivation}.
This is essential in showing that the final frame validates the original logic. Step (3)
can also be utilized to show completeness for the monadic second-order theory of $\omega$
\cite{riba2012model}.

In Section \ref{sec:prelim} we review the basic definitions and provide a Tarski-Vaught-style test for the expansion of the propositional domain to the full powerset to preserve the value of all formulas with propositional quantifiers. In Section \ref{sec:quantifiable}, we show that the saturated canonical general frame of any \pqml\ containing \ttt{Bc} validates the logic. In Section \ref{sec:finite-diversity}, we show how finite diversity allows a quantifiable frame to pass the Tarski-Vaught test. In Section \ref{sec:completeness} we prove our main result. Finally, we conclude in Section \ref{sec:conclusion}.

\section{Preliminaries} \label{sec:prelim}
\begin{definition}
We fix a countably infinite set $\Prop$ of propositional variables and define the language
$\mathcal{L}$ of \pqml s by the following grammar:
\begin{align*}
  \varphi ::= p \mid \lnot\varphi \mid (\varphi \lor \varphi) \mid \Diamond\varphi \mid \exists p\varphi
\end{align*}
where $p \in \Prop$ and $\Diamond$ is the sole modality in this language. The other common
connectives $\land, \to, \leftrightarrow, \Box$ and the universal quantifier $\forall p$
are defined as abbreviations as usual. We also define iterated modalities recursively:
$\Diamond^{0}\varphi = \varphi$, $\Diamond^{n+1}\varphi = \Diamond\Diamond^n\varphi$,
$\Diamond^{\le 0} \varphi = \varphi$, and $\Diamond^{\le n+1}\varphi = \Diamond^{\le
n}\varphi \lor \Diamond^{n+1}\varphi$. $\Box^n$ and $\Box^{\le n}$ are defined dually. Let $\fv(\varphi)$ be the set of free variables of
$\varphi$ and $\mathcal{L}_{qf}$ the quantifier-free fragment of $\mathcal{L}$.
\end{definition}
\begin{definition}
A \emph{Kripke frame} is a pair $\bbf = (W, R)$ where $W$ is non-empty and $R \subseteq
W^2$. A \emph{general frame} is a triple $\bff = (W, R, B)$ where $(W, R)$ is a Kripke
frame and $B \subseteq \wp(W)$ is non-empty and closed under Boolean operations and
$m_\Diamond^\mathbf{F}$, defined by $m_\Diamond^\mathbf{F}(X) = \{w \in W \mid \exists u
\in X, wRu\}$. $R[w] = \{v \in W \mid wRv\}$ and $R[X] = \bigcup_{w \in X} R[w]$. When $B = \wp(W)$ we call $\bff$ \emph{full}.
  
A valuation $v$ for $\bff$ (or $B$) is a function from $\mathsf{Prop}$ to $B$. We define
the semantic value $\tset{\varphi}^{\mathbf{F}}(v)$ of formulas $\varphi$ relative to
valuation $v$ in the general frame $\mathbf{F} = (W, R, B)$ inductively by 
  \begin{align*}
    &\tset{p}^{\mathbf{F}}(v) = v(p) \quad \tset{\lnot\varphi}^{\mathbf{F}}(v) = W \setminus \tset{\varphi}^{\mathbf{F}}(v) \quad \tset{\varphi \lor \psi}^{\mathbf{F}}(v) = \tset{\varphi}^{\mathbf{F}}(v) \cup \tset{\psi}^{\mathbf{F}}(v) \\
    &\tset{\Diamond\varphi}^{\mathbf{F}}(v) = m_\Diamond^\mathbf{F}(\tset{\varphi}^\mathbf{F}(v)) \quad \tset{\exists p \varphi}^{\mathbf{F}}(v) = \bigcup\{\tset{\varphi}^{\mathbf{F}}(v[X/p]) \mid X \in B\}.
  \end{align*}
Here $v[X/p]$ is the function that is identical to $v$ except that $v[X/p](p) = X$. The intended meaning of
$\tset{\varphi}^{\bff}(v)$ is that it is the set of worlds where $\varphi$ is true. A
formula $\varphi$ is \emph{valid} on $\bff$ when $\tset{\varphi}^{\bff}(v) = W$ for all
valuations $v$ for $\bff$.

Valuation and semantics for any Kripke frame $\bbf = (W, R)$ is defined the same as for
$(W, R, \wp(W))$. That is, semantically a Kripke frame is equivalent to the full general
frame based on it.

A \emph{quantifiable frame} is a general frame $\mathbf{F} = (W, R, B)$ such that for any
$v \in B^{\Prop}$ and any $\varphi \in \mathcal{L}$, $\tset{\varphi}^{\mathbf{F}}(v) \in
B$. That is, $B$ is `closed under semantics'.
\end{definition}

\noindent\textbf{Notation:} Given a Kripke frame $(W, R)$ and two general frames
$\mathbf{F} = (W, R, A)$ and $\mathbf{G} = (W, R, B)$ based on it, for any valuation $v
\in (A \cap B)^{\Prop}$, it is clear that if $\tset{\varphi}^\mathbf{F}(v) \not=
\tset{\varphi}^\mathbf{G}(v)$, it is only because of the difference between $A$ and $B$.
Hence, when it is clear which Kripke frame $(W, R)$ is in discussion, we write
$\tset{\varphi}^{(W, R, A)}$ simply as $\tset{\varphi}^A$ or even $\tset{\varphi}$ when
$\varphi$ is quantifier-free. Also, it is routine to show that $\tset{\varphi}^{\bff}(v)$
only depends on $v|_{\fv(\varphi)}$, the restriction of $v$ to $\fv(\varphi)$. Thus, for
any partial function $v$ from $\Prop$ to $B$ such that $dom(v) \supseteq \fv(\varphi)$, we
take $\tset{\varphi}^\bff(v)$ to be the unique element in $\{\tset{\varphi}^\bff(v') \mid
v \subseteq v' \in B^{\Prop}\}$.

\vspace{0.5em}
Now we present the Tarski-Vaught-style test for expanding the propositional domain safely.
\begin{definition}
Given a Kripke frame $(W, R)$ and $\varnothing \not= A, B \subseteq \wp(W)$, $A$ is a
\emph{$\pqml$-invariant subdomain} of $B$ if $A \subseteq B$ and for any $\varphi \in
\mathcal{L}$ and $v \in A^{\Prop}$, $\tset{\varphi}^A(v) = \tset{\varphi}^B(v)$.
\end{definition}
\begin{lemma}\label{lem:Tarski-Vaught}
Given a Kripke frame $(W, R)$ and $\varnothing \not= B \subseteq \wp(W)$, $B$ is a
$\pqml$-invariant subdomain of $\wp(W)$ iff for any $\varphi \in \mathcal{L}$, $p \in
 \fv(\varphi)$, $v \in B^{\mathsf{Prop}}$, $w \in W$, and
$X \in \wp(W)$, if $w \in
\tset{\varphi}^{\wp(W)}(v[X/p])$ then there is $Y \in B$ with $w \in
\tset{\varphi}^{\wp(W)}(v[Y/p])$.
\end{lemma}
\begin{proof}
Left-to-Right: suppose $B$ is a $\pqml$-invariant subdomain of $\wp(W)$ and $w \in
\tset{\varphi}^{\wp(W)}(v[X/p])$. Then $w \in \tset{\exists p \varphi}^{\wp(W)}(v)$. Then
by assumption, $w \in \tset{\exists p \varphi}^B(v)$ and hence $w \in
\bigcup\{\tset{\varphi}^B(v[Y/p]) \mid Y \in B\}$. So there is $Y \in B$ such that $w \in
\tset{\varphi}^B(v[Y/p])$. Finally since $B$ is a $\pqml$-invariant subdomain of $\wp(W)$, and
$v[Y/p]\in B^{\Prop}$, $\tset{\varphi}^B(v[Y/p]) = \tset{\varphi}^{\wp(W)}(v[Y/p])$.
Consequently there is $Y \in B$ such that $w \in \tset{\varphi}^{\wp(W)}(v[Y/p])$. 

Right-to-Left: assume the stated criteria and use induction on $\varphi$. Only the step for $\exists$
is non-trivial, where we need to show that $\bigcup\{\tset{\varphi}^{\wp(W)}(v[X/p]) \mid
X \in \wp(W)\} = \bigcup\{\tset{\varphi}^B(v[Y/p]) \mid Y \in B\}$ with $v \in B^{\Prop}$.
By IH, we only need to show that $\bigcup\{\tset{\varphi}^{\wp(W)}(v[X/p]) \mid X \in
\wp(W)\} = \bigcup\{\tset{\varphi}^\wp(W)(v[Y/p]) \mid Y \in B\}$, and the right-to-left
 inclusion is
trivial. For the other inclusion, take any $w \in
 \bigcup\{\tset{\varphi}^{\wp(W)}(v[X/p])
\mid X \in \wp(W)\}$ and use the criteria.
\end{proof}

Substitutions play an important role in our later proofs. Many authors define
substitution only when it is free to do so, but we need the version that renames the bound
variables when conflicts arise.
\begin{definition}
A \emph{substitution} is a function $\sigma: \Prop \to \mathcal{L}$, and $\sigma^\psi_p$ for any $p \in \Prop$ and $\psi \in \mathcal{L}$ is the substitution that is identical with $\sigma$ except that $\sigma^\psi_p(p) = \psi$. Given a substitution
$\sigma$, we extend it to $\mathcal{L}$ recursively so that that $\sigma(\lnot\varphi) =
\lnot\sigma(\varphi)$, $\sigma(\Diamond\varphi) = \Diamond\sigma(\varphi)$,
$\sigma(\varphi \lor \psi) = \sigma(\varphi) \lor \sigma(\psi)$, and $\sigma(\exists p
\varphi) = \exists q \sigma_p^q(\varphi)$ where $q = p$ if $p \not\in \bigcup_{r \in
\fv(\exists p \varphi)}{\fv(\sigma(r))}$ and otherwise $q$ is the first variable not used
in $\exists p \varphi$ and any $\sigma(r)$ for $r \in \fv(\exists p\varphi)$.  
\end{definition}
\noindent Let $\iota$
be the identity substitution. Then, using the above definition, $\iota_p^\psi(\varphi)$ is the result of substituting $\psi$ for $p$ in $\varphi$
with the necessary renamings of bound variables. In particular, $\iota(\varphi) =
\varphi$. Then, the standard substitution lemma connecting syntactic and semantic
substitution is:
\begin{lemma}
On any general frame $\bff = (W, R, B)$, valuation $v$ for $\bff$, and substitution
$\sigma$, define valuation $\sigma\star v: p \mapsto \tset{\sigma(p)}^\bff(v)$. Then for
any $\varphi \in \mathcal{L}$, $\tset{\varphi}^\bff(\sigma \star v) =
\tset{\sigma(\varphi)}^\bff(v)$.
\end{lemma}

Finally, we introduce logic. For the convenience of certain proofs, we opted for
$\Diamond$ as the primitive modality. For this reason, the \ttt{Dual} axiom is necessary.
\begin{definition}
A \emph{normal propositionally quantified modal logic} (\npqml) is a set $\Lambda \subseteq
\mathcal{L}$ satisfying the following conditions:
  \begin{itemize}
\item (Taut) all propositional tautologies are in $\Lambda$;
\item axiom $\mathtt{K} = \Box(p \to q) \to (\Box p \to \Box q)$ and $\mathtt{Dual} =
\Diamond p \leftrightarrow \lnot\Box\lnot p$ are in $\Lambda$;
\item (EI) all instances of $\iota_p^{\psi}(\varphi) \to \exists p \varphi$ are in
$\Lambda$;
\item (Nec) if $\varphi \in \Lambda$, then $\Box\varphi \in \Lambda$;
\item (MP) if $\varphi, (\varphi \to \psi) \in \Lambda$, then $\psi \in \Lambda$;
\item (EE) if $\varphi \to \psi \in \Lambda$ with $p \not\in \fv(\psi)$, then
$\exists p\varphi \to \psi \in \Lambda$.
  \end{itemize}
For any axioms or axiom schemes $A_1, A_2, \dots, A_n$, we write $\mathsf{K}_\Pi A_1
A_2\dots A_n$ for the smallest \npqml\ containing all (instances) of all $A_i$'s.
\end{definition}
\begin{fact}
For any $\varphi$ and $\psi$ obtained by renaming some bound variables in $\varphi$, $(\varphi \leftrightarrow \psi) \in \Lambda$ for any \npqml\ $\Lambda$.
\end{fact}
Recall that the famous \emph{Barcan scheme} \ttt{Bc} is $\Diamond\exists p \varphi \to
\exists p\Diamond\varphi$. 
\begin{fact}
For any class $\mathsf{C}$ of quantifiable frames, the set $\Lambda$ of formulas valid
on each member of $\mathsf{C}$ is a \npqml\ containing all instances of \ttt{Bc}.
\end{fact}
\begin{definition}
For any set $\Theta \subseteq \mathcal{L}$, let $\mathsf{KFr}(\Theta)$ be the class of Kripke frames validating all formulas in $\Theta$ and let $\Theta\pi+$ be the set of formulas valid on all members of $\mathsf{KFr}(\Theta)$.
\end{definition}

\section{General completeness for quantifiable frames} \label{sec:quantifiable}

This section introduces saturated canonical general frames for \npqml s containing all
instances of the Barcan scheme and shows that they validate the original logic and are
automatically quantifiable. A consequence is the following:
\begin{theorem}\label{thm:completeness-quantifiable-frames} Any \npqml\ $\Lambda \supseteq
\mathsf{K}_\Pi\mathtt{Bc}$ is strongly complete w.r.t. the class of quantifiable frames
validating $\Lambda$.
\end{theorem}
Fix a \npqml\phantom{ }  $\Lambda \supseteq \mathsf{K}_\Pi\mathtt{Bc}$. To construct the saturated
canonical general frame for $\Lambda$, we extend $\Prop$ to $\Prop^+$ with countably
infinitely many new variables and obtain the extended language $\mathcal{L}^+$. Semantics
and logics for $\mathcal{L}^+$ are defined completely analogously. Let $\Lambda^+$ be the
smallest \npqml\ in $\mathcal{L}^+$ extending $\Lambda$.
\begin{definition}
A set $\Gamma \subseteq \mathcal{L}$ is \emph{$\Lambda$-consistent} if there is no finite
$A \subseteq \Gamma$ s.t. $\lnot(\bigwedge A) \in \Lambda$. A \emph{maximally
$\Lambda$-consistent set} ($\Lambda$-MCS) is a $\Lambda$-consistent set s.t. all of its proper
extensions in $\mathcal{L}$ are not $\Lambda$-consistent. $\Lambda^+$-consistency and
$\Lambda^+$-MCSs are defined in the same way using $\Lambda^+$ and $\mathcal{L}^+$. 

A $\Lambda^+$-MCS $\Gamma$ is \emph{saturated} if for any $\exists p\varphi \in \Gamma$,
there is $q \in \Prop^+$ not occurring in $\exists p\varphi$ s.t. $\iota_p^q(\varphi) \in
\Gamma$.
\end{definition}
Now define the \emph{saturated canonical general frame} $\bff_\Lambda = (W, R, B)$ where 
\begin{itemize}
\item $W$ is the set of all saturated $\Lambda^+$-MCSs, 
\item $w R v$ iff for all $\varphi \in v$, $\Diamond\varphi \in w$, 
\item $B = \{[\varphi] \mid \varphi \in \mathcal{L}^+\}$ where $[\varphi] = \{w \in W \mid
\varphi \in w\}$.
\end{itemize}
The following two lemmas are completely analogous to their first-order modal logic
counterparts. The proof of the first extension lemma is almost identical for example to
the proof of Theorem 14.1 of \cite{Cresswell1996-CREANI-3}. The full power of (EI) is not
used, and all we need is (EE) and that \npqml s can prove equivalences between formulas
that differ only by renaming of bound variables. The second existence lemma can also be
proved by repeating the steps in the proof of the existence lemma for first-order modal
logic with \ttt{Bc}. For an example, see the proof of Theorem 14.2 of
\cite{Cresswell1996-CREANI-3}.
\begin{lemma}
Any $\Lambda$-consistent set of $\mathcal{L}$ formulas $\Gamma \subseteq \mathcal{L}$ can be extended to a saturated $\Lambda^+$-MCS
$\Gamma^+$.
\end{lemma}
\begin{lemma}\label{lem:existence-general} For any $w \in W$, if $\Diamond\varphi \in w$,
then there is $u \in R[w]$ s.t. $\varphi \in u$.
\end{lemma}

The truth lemma is replaced by a more general statement for all valuations arising from
substitutions. This is a common idea in algebraic semantics.
\begin{lemma}\label{lem:truth-substitution} For any substitution $\sigma$ for
$\mathcal{L}^+$, define its associated valuation $[\sigma]: p \mapsto [\sigma(p)]$. Then
for any $\varphi \in \mathcal{L}^+$ and all $\sigma$,
$\tset{\varphi}^{\bff_\Lambda}([\sigma]) = [\sigma(\varphi)]$.
\end{lemma}
\begin{proof}
By induction on $\varphi$. The cases for variables and negation go by 
  \begin{align*}
    &\tset{p}^{\bff_\Lambda}([\sigma]) = [\sigma](p) = [\sigma(p)]. \\
    &\tset{\lnot \varphi}^{\bff_\Lambda}([\sigma]) = W \setminus \tset{\varphi}^{\bff_\Lambda}([\sigma]) = W \setminus [\sigma(\varphi)] = [\lnot\sigma(\varphi)] = [\sigma(\lnot\varphi)].
  \end{align*}
The case for disjunction is similar. For the modal case, note that Lemma
\ref{lem:existence-general} implies that for any $\varphi$, $[\Diamond\varphi] =
m^{\bff}_\Diamond([\varphi])$, so this is again easy.

For the quantifier case, recall that $\sigma(\exists p \varphi) = \exists q
\sigma_p^q(\varphi)$ for a suitable $q$. Now observe that due to saturation and (EI),
  \begin{align*}
    [\exists q \sigma_p^q(\varphi)] \subseteq \bigcup\{[\iota_q^r\sigma_p^q(\varphi)] \mid r \in \Prop^+\} \subseteq \bigcup\{[\iota_q^\psi\sigma_p^q(\varphi)] \mid \psi \in \mathcal{L}^+\} \subseteq [\exists q \sigma_p^q(\varphi)]. 
  \end{align*}
Next, observe that given how $q$ is chosen when performing $\sigma(\exists p \varphi) =
\exists q \sigma_p^q(\varphi)$, $\iota_q^\psi\sigma_p^q(\varphi)$ and
$\sigma_p^\psi(\varphi)$ differ only by renaming of bound variables and thus are logically
equivalent. Hence, $[\iota_q^\psi\sigma_p^q(\varphi)] = [\sigma_p^\psi(\varphi)]$. Also, $[\sigma_p^\psi] = [\sigma][[\psi]/p]$. Thus, with IH and recalling that
$B = \{[\psi] \mid \psi \in \mathcal{L}^+\}$, 
  \begin{align*}
    [\sigma(\exists p \varphi)] &= \bigcup\{[\sigma_p^\psi(\varphi)] \mid \psi \in \mathcal{L}^+\} = \bigcup\{\tset{\varphi}^{\bff_\Lambda}( [\sigma_p^\psi] ) \mid \psi \in \mathcal{L}^+\} \\
    &= \bigcup\{\tset{\varphi}^{\bff_\Lambda}([\sigma][[\psi]/p]) \mid \psi \in \mathcal{L}^+\} = \bigcup\{\tset{\varphi}^{\bff_\Lambda}([\sigma][X/p]) \mid X \in B\} \\
    &= \tset{\exists p \varphi}^{\bff_\Lambda}([\sigma]). 
  \end{align*}

  \vspace{-2.2em}
  \phantom{xxx}
\end{proof}

Now we are ready to prove Theorem \ref{thm:completeness-quantifiable-frames}. First we show that $\bff_{\Lambda}$ is indeed a quantifiable frame, for any valuation $v\in B^\Prop$ for $\bff_\Lambda$, given how $B$ is defined, consider the substitution $\sigma_v:p\mapsto \varphi_p$ where $v(p) = [\varphi_p]$. Then $v = [\sigma_v]$. It follows that $\bff_\Lambda$ is quantifiable by Lemma \ref{lem:truth-substitution}. Next, $\bff_{\Lambda}\models \Lambda$ as for $\varphi \in \Lambda \subseteq \Lambda^+$ and arbitrary $v\in B^\Prop$, $[\varphi]^{\bff_\Lambda}(v) = [\varphi]^{\bff_\Lambda}(\sigma_v) = [\sigma_v(\varphi)] = W$ since $\sigma(\varphi)$ is also in $\Lambda^+$. Finally, taking $\sigma$ as $\iota$ in Lemma \ref{lem:truth-substitution}, $\tset{\varphi}^{\bff_\Lambda}([\iota]) = [\varphi]$, which means under valuation $[\iota]$, each $\Lambda^+$-MCS is satisfied by itself on a general frame $\bff_\Lambda$ that validates $\Lambda$. Hence Theorem \ref{thm:completeness-quantifiable-frames} follows.

\section{From finite diversity to {\normalfont\textsc{pqml}}-invariant subdomain} \label{sec:finite-diversity}
We first introduce the concepts of duplicates and diversity. 
\begin{definition}
Given a Kripke frame $\bbf = (W, R)$, we say that worlds $w, u \in W$ are
\emph{duplicates} if the permutation $(wu)$ of $W$ that exchanges $w$ and $u$ is an
automorphism of $\bbf$. Let $\Delta$ be this relation of being duplicates ($\bbf$'s
duplication relation), which clearly is an equivalence relation on $W$, and then let
$W/\Delta$ be the set of $\Delta$'s equivalence classes. The \emph{diversity} of $\bbf$ is the
cardinality of $W/\Delta$. The \emph{diversity} of a Kripke frame class is the supremum of
the diversity of all point-generated subframes of the frames in that class (if exists). 
\end{definition}

Intuitively, duplicate classes are `positions' a world could be in, and the diversity of a Kripke frame counts the number of positions in that frame. We use some examples to illustrate the concept of diversity. 
\begin{example}
While cyclic frames of the form $(\{0, 1, \dots, n-1\}, \{(i, i+1 \mathrel{mod} n) \mid i
= 0 \dots n-1\})$ are highly symmetric, no two distinct worlds are duplicates of each
other, as switching them and only them is not an automorphism.
\end{example}
\begin{example}
    The diversity of $\mathsf{KFr}(\mathtt{D45})$ is 2, and the diversity of $\mathsf{KFr}(\mathtt{5})$ is 3. It is well known that a point-generated Kripke frame validating $\mathtt{D45}$ is either a clique $(W, W \times W)$ or a point looking at a clique $(\{r\} \cup W, (\{r\} \cup W) \times W)$. Clearly, Kripke frames of the later form has exactly two duplicate classes as every $x, y \in W$ are duplicates of each other.

    It is also well known that the most non-trivial kind of point-generated Kripke frames validating $\mathtt{5}$ are of the form $(\{r\} \cup W, \{r\} \times U \cup W \times W)$ where $r \not \in W$ and $\varnothing \not= U \subseteq W$. Then, there are three duplicate classes: $\{r\}, U, W \setminus U$. 
\end{example}
\begin{example}
    We give a Sahlqvist definable frame class of diversity $4$. Consider Sahlqvist formulas $\varphi_1 = \Box(\Diamond p \to \Box\Diamond p)$ and $\varphi_2 = \Diamond\Diamond p \to \Box\Diamond p$. Let $\mathbb{F} = (W, R)$ be a Kripke frame that validates the two formulas and is point-generated from $r \in W$. By $\varphi_2$, for any $x, y \in R[r]$, $R[x] = R[y]$. By $\varphi_1$, all worlds $w \in W \setminus \{r\}$ locally validates $\mathtt{5}$: $(W, R), w \models \forall p (\Diamond p \to \Box\Diamond p)$. Now we discuss several cases:
    \begin{itemize}
        \item If $R[r] = \varnothing$, then $(W, R) = (\{r\}, \varnothing)$ with diversity $1$.
        \item If $r \in R[r]$, then $(W, R)$ is a universally connected clique, again with diversity $1$.
        \item If we are not in the above two cases, and there is $u \in R[r]$ that is reflexive, let $A$ be the set of reflexive worlds in $R[r]$, $B$ be the irreflexive worlds in $R[r]$, and $C$ be $R[u]$. Then observe that $W$ is $\{r\} \cup B \cup C$ where $r \not \in B \cup C$, $B \cap C = \varnothing$, and $A \subseteq C$, and $R$ is $\{r\} \times (B \cup A) \cup (B \cup C) \times C$. Clearly, such a frame has at most $4$ duplicate classes: $\{r\}$, $B$ (could be empty), $A$, and $C \setminus A$ (could be empty). 
        \item If we are not in the above three cases, then $R[r]$ is non-empty (call it $A$ and let $u$ be a member of it), and every world in $R[r]$ is irreflexive. If $R[u]$ is empty, then $(W, R)$ is $(\{r\} \cup A, \{r\} \times A)$. If $R[u]$ is non-empty, call it $B$, and let $C$ be $R[v]$ for any $v \in B$. The choice of $v$ is irrelevant due to the $\mathtt{5}$ axiom. Thus, $(W, R)$ is $(\{r\} \cup A \cup C, \{r\} \times A \cup A \times B \cup C \times C)$ where $r \not \in A \cup C$, $A \cap C = \varnothing$, and $B \subseteq C$. Again, $(W, R)$ has at most $4$ duplicate classes: $\{r\}$, $A$, $B$ and $C\setminus B$.
    \end{itemize}
\end{example}
\begin{example}
Finite frames have only finite diversity. Thus, the axioms $\mathtt{Alt}^n
= \bigwedge_{i < n+1}p_i \to \bigvee_{i < j < n+1} \Diamond(p_i \land p_j)$ and $\mathtt{Trs}^m = \Diamond^{\le m}p \to \Diamond^{\le m+1}p$ together define Kripke
frame classes of diversity at most $n^{m+1}$. 
\end{example}

The following lemma collects some easy but very useful properties of the duplicate classes
and how they interact with $R$ and $m_\Diamond$.
\begin{lemma}\label{lem:diversity-analysis} For any Kripke frame $\bbf = (W, R)$ and its
duplication relation $\Delta$,
  \begin{itemize}
\item for any $D_1 \not= D_2 \in W/\Delta$, there is $w \in D_1$ and $u \in D_2$ s.t.
$wRu$ iff for all $w \in D_1$ and $u \in D_2$, $wRu$; 
\item for any $D \in W/\Delta$, the only possible configurations for $R|_D$ are: $D^2$,
$\varnothing$, and when $|D| \ge 2$, $D^2 \setminus id_D$ and $id_D$ ($id_D$ is the
identity relation on $D$).
  \end{itemize}
For convenience we define the binary relation $R_\Delta$ on $W/\Delta$ s.t. $D_1 R_\Delta
D_2$ iff there is $w \in D_1$ and $u \in D_2$ s.t. $wRu$. Here we allow $D_1 = D_2$.

Now we discuss the possible ways $m_\Diamond^{\bff}(X) \cap D$ is determined.
  \begin{itemize}
\item In case $R|_D = D^2$ or $R|_D = \varnothing$, clearly $m_\Diamond^{\bff}(X) \cap D$
is either $D$ or $\varnothing$, and it is $D$ iff there is $D' \in R_\Delta[D]$ s.t. $|X
\cap D'| \ge 1$.
\item In case $R|_D = D \setminus id_D$ with $|D| \ge 2$, if there is $D' \in
R_{\Delta}[D] \setminus \{D\}$ s.t. $|X \cap D'| \ge 1$, then $m_\Diamond^{\bff}(X) \cap D
= D$, otherwise, 
    \begin{itemize}
\item if $|X \cap D| \ge 2$, then also $m_\Diamond^{\bff}(X) \cap D = D$,  
\item if $|X \cap D| = 1$, then $m_\Diamond^{\bff}(X) \cap D = D \setminus X$, and
\item if $|X \cap D| = 0$, then $m_\Diamond^{\bff}(X) \cap D = \varnothing$.
    \end{itemize}
\item In case $R|_D = id_D$ with $|D| \ge 2$, if there is $D' \in R_{\Delta}[D] \setminus
\{D\}$ s.t. $|X \cap D'| \ge 1$, then $m_\Diamond^{\bff}(X) \cap D = D$, and otherwise,
$m_\Diamond^{\bff}(X) \cap D = X$.
  \end{itemize}
\end{lemma}

What follows is the core of the proof of our main result. We want to show that if the
underlying Kripke frame of the quantifiable frame $(W, R, B)$ has finite diversity and $B$
contains all singletons and duplicate classes, then $B$ is a $\pqml$-invariant subdomain of
$\wp(W)$. The key idea is that whenever $\varphi$ is true at $w$ where at most one $p \in
\fv(\varphi)$ is evaluated to a set $X \subseteq W$ that is not necessarily in $B$, we can
always swap the valuation of $p$ to a $Y \in B$ while keeping $\varphi$ true at $w$. For
this to be true, we must establish that $\varphi$'s truth at $w$ is insensitive to certain
changes in the valuation of $p$. For monadic second-order logic, this can be done with
EF-game, but for modal logic, we cannot only focus on how $\varphi$'s truth at $w$ is
insensitive to change since modality requires us to also consider the truth of $\varphi$
at other worlds. We must take a more global perspective and strive to show that
$\varphi$'s `meaning' is insensitive to certain changes. In the end, we arrive at a
qualified quantifier-elimination: when restricted to a duplicate class $D$ and relative to
a valuation $v$, the `meaning' of $\varphi$ can be written as a Boolean formula
$f_\varphi(v, D)$ using variables in $\fv(\varphi)$, and when valuations $u$ and $v$ are
close enough, $f_\varphi(v, D) = f_\varphi(u, D)$.

\begin{definition}
For any finite $\sfp \subseteq \Prop$, we write $\langle \sfp\rangle$ for the Boolean
language generated from $\sfp$. $at(\langle \sfp \rangle)$ is the finite set of all
formulas $l_1 \land \dots \land l_k$ where each $l_i$ is either $p_i$ or $\lnot p_i$ and
$p_1, \dots, p_k$ list all elements in $\sfp$. These formulas correspond to the atoms in
the Lindenbaum algebra of $\langle \sfp \rangle$.
\end{definition}

\begin{definition}
Given a Kripke frame $(W, R)$ with duplicate relation $\Delta$, for any finite $\sfp
\subseteq \Prop$, $u,v \in \wp(W)^{\sfp}$, and $n \in \mathbb{N}$, $u \approx_n v$ if for
all $D \in W/\Delta$ and $\zeta \in at(\langle \sfp \rangle)$, $|\tset{\zeta}(u) \cap D| =
|\tset{\zeta}(v) \cap D|$ or both $|\tset{\zeta}(u) \cap D|, |\tset{\zeta}(v) \cap D| \ge
2^n$.
\end{definition}
\begin{lemma}\label{lem:approx-extension} Let $u, v \in \wp(W)^{\sfp}$ for some finite $\sfp \subseteq \Prop$. 
\begin{enumerate}
    \item If $u \approx_n v$, then not only for $\zeta \in at(\langle \sfp
    \rangle)$, for any $\beta \in \langle \sfp \rangle$ and $D \in W/\Delta$,
    $|\tset{\beta}(u) \cap D| = |\tset{\beta}(v) \cap D|$ or both $|\tset{\beta}(u) \cap D|,
    |\tset{\beta}(v) \cap D| \ge 2^n$. 
    \item If $u \approx_{n} v$ and $p \not\in \sfp$, for
any $X \in \wp(W)$ there is $Y \in \wp(W)$ s.t. $u[X/p] \approx_{n-1} v[Y/p]$ (since $p
\not\in dom(u)$, $u[X/p] = u \cup \{(p, X)\}$).
\end{enumerate}
\end{lemma}
\begin{proof}
    The first part is easy since every $\tset{\beta}(u)$ is the union of some $\tset{\zeta}(u)$'s where $\zeta\in at(\langle \sfp
    \rangle)$.
The second part: for each $\zeta \in at(\langle \sfp \rangle)$ and $D \in W/\Delta$,
choose a set $Y_{\zeta, D} \subseteq \tset{\zeta}(v) \cap D$ such that:
\begin{itemize}
  \item if $|(\tset{\zeta}(u) \cap D) \cap X| < 2^{n-1}$, then $|Y_{\zeta, D}| = |(\tset{\zeta}(u) \cap D) \cap X|$;
\item if $|(\tset{\zeta}(u) \cap D) \setminus X| < 2^{n-1}$, then $|(\tset{\zeta}(v) \cap
D) \setminus Y_{\zeta, D}| = |(\tset{\zeta}(u) \cap D) \setminus X|$;
\item if both $|(\tset{\zeta}(u) \cap D) \cap X|, |(\tset{\zeta}(u) \cap D) \setminus X|
\ge 2^{n-1}$, then $|Y_{\zeta, D}| = 2^{n-1}$.
\end{itemize}
Given that $u \approx_n v$, the above conditions can be satisfied. Then either
$|(\tset{\zeta}(u) \cap D) \cap X| = |Y_{\zeta}, D|$ or both $|(\tset{\zeta}(u) \cap D)
\cap X|, |Y_{\zeta}, D| \ge 2^{n-1}$, and the same goes for $|(\tset{\zeta}(v) \cap
D) \setminus Y_{\zeta, D}|$ and $|(\tset{\zeta}(u) \cap D) \setminus X|$. Then, with $Y =
\bigcup_{\zeta \in at(\langle \sfp \rangle), D \in W/\Delta}Y_{\zeta, D}$, $u[X/p]
\approx_{n-1} v[Y/p]$.
\end{proof}
Let $qd(\varphi)$ be the quantifier depth of $\varphi$.
\begin{lemma}\label{lem:boolean-breakdown}
Given a Kripke frame $(W, R)$ with duplicate relation $\Delta$, for each $\varphi \in
\mathcal{L}$, there is a function $f_\varphi : (\wp(W)^{\fv(\varphi)} \times W/\Delta) \to
\langle \fv(\varphi) \rangle$ such that 
  \begin{itemize}
\item for any $v \in \wp(W)^{\fv(\varphi)}$ and $D \in W/\Delta$, $\tset{\varphi}(v) \cap
D = \tset{f_\varphi(v, D)}(v) \cap D$;
\item for any $u \approx_{qd(\varphi) + 1} v \in \wp(W)^{\fv(\varphi)}$, $f_\varphi(u, D)
= f_{\varphi}(v, D)$ for all $D \in W/\Delta$. 
  \end{itemize}
In particular, recursively define $f$ as follows and it will witness the lemma: for the
base and Boolean cases: $f_p(v, D) = p$, $f_{\lnot\varphi}(v, D) = \lnot f_{\varphi}(v,
D)$, $f_{\varphi \lor \psi}(v, D) = f_{\varphi}(v|_{\fv(\varphi)}, D) \lor
f_{\psi}(v|_{\fv(\psi)}, D)$. For the modal case, we copy the analysis in Lemma
\ref{lem:diversity-analysis} with $X = \tset{f_\varphi(v, D)}(v)$:
 \begin{itemize}
\item In case $R|_D = D^2$ or $\varnothing$, if there is $D' \in R_\Delta[D]$ s.t. $|X
\cap D'| \ge 1$ then $f_{\Diamond\varphi}(v, D) = \top$, otherwise $f_{\Diamond\varphi}(v,
D) = \bot$.
\item In case $R|_D = D \setminus id_D$ with $|D| \ge 2$, if there is $D' \in
R_{\Delta}[D] \setminus \{D\}$ s.t. $|X \cap D'| \ge 1$, then $f_{\Diamond\varphi}(v, D) =
\top$, otherwise, 
    \begin{itemize}
\item if $|X \cap D| \ge 2$, then also $f_{\Diamond\varphi}(v, D) = \top$,  
\item if $|X \cap D| = 1$, then $f_{\Diamond\varphi}(v, D) = \lnot f_{\varphi}(v, D)$, and
\item if $|X \cap D| = 0$, then $f_{\Diamond\varphi}(v, D) = \bot$.
    \end{itemize}
\item In case $R|_D = id_D$ with $|D| \ge 2$, if there is $D' \in R_{\Delta}[D] \setminus
\{D\}$ s.t. $|X \cap D'| \ge 1$, then $f_{\Diamond\varphi}(v, D) = \top$, and otherwise,
$f_{\Diamond\varphi}(v, D) = f_{\varphi}(v, D)$.
  \end{itemize}
\noindent For the quantifier case let $f_{\exists p \varphi}(v, D)$ be the following:
  \begin{align*}
    \bigvee\{\zeta \in at(\langle \fv(\exists p \varphi) \rangle) \mid \tset{\zeta}(v) \cap \bigcup_{X \in \wp(W)}\tset{f_\varphi(v[X/p], D)}(v[X/p]) \cap D \not= \varnothing\}.
  \end{align*}
\end{lemma}
\begin{proof}
We first show by induction that whenever $u \approx_{qd(\varphi)+1} v$, $f_\varphi(u, D) =
f_\varphi(v, D)$. The base case and the inductive steps for Boolean connectives are
trivial. For the modal case, suppose $u \approx_{qd(\Diamond\varphi)+1} v$. Then $u
\approx_{qd(\varphi) + 1} v$. Using IH, let $\beta = f_\varphi(u, D) = f_\varphi(v, D)$
and let $X = \tset{\beta}(u), X' = \tset{\beta}(v)$. Now at least $u \approx_{1} v$, so
for any $E \in W/\Delta$, either $|X \cap E| = |X' \cap E|$ or both $|X \cap E|$ and $|X'
\cap E| \ge 2$. This means in the case analysis defining $f_{\Diamond\varphi}(u, D)$ and
$f_{\Diamond\varphi}(v, D)$, the same case must be active, and $f_{\Diamond\varphi}(u, D)
= f_{\Diamond\varphi}(v, D)$.

For the quantifier case, suppose $u \approx_{qd(\exists p \varphi) + 1} v$ with $\sfp =
\fv(\exists p \varphi)$ and $u, v \in \wp(W)^{\sfp}$. Then $u \approx_{qd(\varphi) + 2}
v$. Now pick any $\zeta \in at(\langle \sfp \rangle)$ and suppose $\zeta$ is a disjunct of
$f_{\exists p \varphi}(u, D)$. Then there is $X \in \wp(W)$ s.t. $\tset{\zeta}(u) \cap
\tset{f_\varphi(u[X/p], D)}(u[X/p]) \cap D \not= \varnothing$. Since $p \not \in \sfp$,
$\tset{\zeta}(u) = \tset{\zeta}(u[X/p])$. So we have $\tset{\zeta \land f_\varphi(u[X/p],
D)}(u[X/p]) \cap D \not= \varnothing$. Now by Lemma \ref{lem:approx-extension}, there is
$Y \in \wp(W)$ s.t. $u[X/p] \approx_{qd(\varphi)+1} v[Y/p]$. So with IH, we can let $\beta
= \zeta \land f_\varphi(u[X/p], D) = \zeta \land f_\varphi(v[Y/p], D)$ and now
$\tset{\beta}(u[X/p]) \cap D \not=\varnothing$. Then $\tset{\beta}(v[Y/p]) \cap D \not=
\varnothing$. This means $\zeta$ is also a disjunct of $f_{\exists p \varphi}(v, D)$. The
above argument can be reversed, so $f_{\exists p\varphi}(u, D)$ and $f_{\exists p
\varphi}(v, D)$ have the same disjuncts and thus are the same formula.
  
Now we show that $\tset{\varphi}(v) \cap D = \tset{f_{\varphi}(v, D)}(v) \cap D$. Again
this is by induction and the non-quantifier cases are easy. For easy notation, let $\sfp =
\fv(\exists p \varphi)$ and $\beta_X = f_\varphi(v[X/p], D)$. By IH, $\tset{\exists p
\varphi}(v) \cap D = \bigcup_{X \in \wp(W)} \tset{\beta_X}(v[X/p]) \cap D$. Given the
definition of $f_{\exists p \varphi}(v, D)$ and that $\mathcal{C} := \{\tset{\zeta}(v)
\cap D \mid \zeta \in at(\langle \sfp \rangle)\}$ forms a partition of $D$, all we need to
show is that $\tset{\exists p \varphi}(v) \cap D$ is a union of cells in $\mathcal{C}$.
For this, it suffices to show that for any $\zeta \in at(\langle \exists p
\varphi\rangle)$ and $w_1, w_2 \in \tset{\zeta}(v) \cap D$, if $w_1 \in \bigcup_{X \in
\wp(W)} \tset{\beta_X}(v[X/p]) \cap D$ then $w_2$ is also in $\bigcup_{X \in \wp(W)}
\tset{\beta_X}(v[X/p]) \cap D$. So suppose that $w_1, w_2 \in \tset{\zeta}(v) \cap D$ for
some $\zeta \in at(\langle p \rangle)$ and there is $X \in \wp(W)$ s.t. $w_1 \in
\tset{\beta_X}(v[X/p]) \cap D$. Recall that $(w_1w_2)$ is the permutation of $W$ that
exchanges $w_1$ and $w_2$. Let $Y = (w_1w_2)[X]$. Since $w_1, w_2$ are both in
$\tset{\zeta}(v)$ and $\zeta$ is in $at(\langle \sfp \rangle)$, for any $q \in \sfp$,
$v(q) = (w_1w_2)[v(q)]$ as $w_1 \in v(q)$ iff $w_2 \in v(q)$. Recall also that $w_1, w_2$
are in the same duplication class $D$, so for any $D' \in W/\Delta$, $w_1 \in D'$ iff $w_2
\in D'$. From all these, it is clear that $v[X/p] \approx_{qd(\varphi) + 1} v[Y/p]$, since
$(w_1w_2)[\cdot]$ commutes with all Boolean operations and thus for any $\gamma \in
\langle \sfp \cup \{p\} \rangle$ and $D' \in W/\Delta$, $\tset{\gamma}(v[Y/p]) \cap D' =
(w_1w_2)[\tset{\gamma}(v[X/p]) \cap D']$, meaning also that they are of the same
cardinality. Then $\beta_X = \beta_Y$. Recall that $w_1 \in \tset{\beta_X}(v[X/p]) \cap
D$. Apply $(w_1w_2)$ to both sides and we have $w_2 \in \tset{\beta_Y}(v[Y/p]) \cap D$.
\end{proof}

\begin{theorem}\label{thm:ele-subdomain} For any general frame $(W, R, B)$ with
$(W, R)$ having finite diversity, $B$ is a $\pqml$-invariant subdomain of $\wp(W)$ if for any $w \in W$ and $D \in W/\Delta$, $\{w\}, D \in B$. (Only $B$'s closure under Boolean operations is used.)
\end{theorem}

\begin{proof}
It suffices to show that for any $w \in W$, $\varphi \in \mathcal{L}$ with $\sfp =
\fv(\exists p \varphi)$ and $n = qd(\varphi) + 1$, $v \in B^\sfp$, and $X \in \wp(W)$,
there is $Y \in B$ s.t. $v[X/p] \approx_n v[Y/p]$ and $w \in X$ iff $w \in Y$, since if
so, then by Lemma \ref{lem:boolean-breakdown}, $w\in \tset{\varphi}(v[X/p])\cap D$ iff
$w\in \tset{f_\varphi(v[X/p],D)}(v[X/p])\cap D$ iff $w\in
\tset{f_\varphi(v[Y/p],D)}(v[X/p])\cap D$ iff $w\in \tset{f_\varphi(v[Y/p],D)}(v[Y/p])\cap
D$ (recall that $f_\varphi(v[Y/p], D)$ is a Boolean formula) iff $w\in
\tset{\varphi}(v[Y/p])\cap D$ and thereby by Lemma \ref{lem:Tarski-Vaught} we are done.

Notice that
by assumption we have a finite partition $\mathcal{C} = \{\tset{\zeta}(v) \cap D \mid
\zeta \in at(\langle \sfp \rangle), D \in W/\Delta\}$ of $W$ and each $C \in \mathcal{C}$
is in $B$ by assumption. For $Y \in B$, it is enough to make sure that for all $C \in
\mathcal{C}$, $C \cap Y$ or $C \setminus Y$ is finite. For $v[X/p] \approx_n v[Y/p]$, it
is enough to make sure that for all $C \in \mathcal{C}$, $|C \cap Y| = |C \cap X|$ or $|C
\cap Y|, |C \cap X| \ge 2^n$, and $|C \setminus Y| = |C \setminus X|$ or $|C \setminus Y|,
|C \setminus X| \ge 2^n$. Thus, for each $C \in \mathcal{C}$, let $Y_C = C \cap X$ if
either $C \cap X$ or $C \setminus X$ is finite, and otherwise when both are infinite, if
$w \in C \cap X$, let $Y_C$ be $C \setminus Z$ for some $Z \subseteq C \setminus X$ with
$|Z| = 2^n$, and otherwise let $Y_C$ be some subset of $C \cap X$ with $|Y_C| = 2^n$.
Clearly $Y = \bigcup_{C \in \mathcal{C}} Y_C$ satisfies the requirements.
\end{proof}

\section{General completeness with finite-diversity}\label{sec:completeness}
We first define the extra axioms needed. 
\begin{definition}
  \begin{itemize}
\item $\mathtt{Q}^n(\varphi) = \Diamond^{\le n}\varphi \land \forall p (\Box^{\le
n}(\varphi \to p) \lor \Box^{\le n}(\varphi \to \lnot p))$ where $p$ is the first variable
not in $\fv(\varphi)$;
\item $\mathtt{At}^n = \forall q (\Diamond^{\le n} q \to \exists p (\mathtt{Q}^{\le n}(p)
\land \Box^{\le n}(p \to q)))$;
\item $\mathtt{R}^{n} = \forall p \exists q (\Box^{\le n}(p \to \Box q) \land \forall
r (\Box^{\le n}(p \to \Box r) \to \Box^{\le n}(q \to r))$.
  \end{itemize}
\end{definition}
\noindent The iterated modalities $\Diamond^{\le n}$ and $\Box^{\le n}$ used here are
meant to capture the global modality (at least at the root of a point-generated model).
Then, $\mathtt{Q}^n(\varphi)$ states that the proposition $\varphi$ expresses is a
maximally specific possible proposition, a non-bottom proposition that settles the truth
or falsity of all propositions. Semantically, assuming that enough propositions exist to
distinguish possible worlds, such maximally specific possible propositions must be
singletons, true at exactly one possible world. Even without assuming that worlds can be
distinguished by propositions, maximally specific possible propositions are atoms in the
algebra of all propositions, and can themselves serve as worlds. Thus, they are often
called world propositions. Using $\mathtt{Q}^n$, $\mathtt{At^n}$ is the atomicity
principle stating that whatever is possible is entailed by a maximally specific possible
proposition. In other words, the algebra of all propositions is atomic. Such an atomicity
principle features prominently in many works
\cite{fine1970propositional,ten2006expressivity,Holliday2019note,Ding2018-DINOTL,fritz2023foundations}
with both technical and philosophical significance. The use of the nested modality
$\Box^{\le n}$ is needed as we are not assuming strong modal axioms and we need to
simulate a global modality. The formula $\mathtt{R^n}$ is meant to capture the fact that
$\wp(W)$ is trivially closed under taking successor sets: if $X \in \wp(W)$, then $R[X]
\in \wp(W)$. Again, using $\Box^{\le n}$ as a substitute for the global modality,
$\mathtt{R}^n$ says that for any proposition $p$, there is a proposition $q$ such that
whenever $p$ is true, `necessarily $q$' is true, and $q$ is the strongest proposition with
this property. The validity of these two formulas over Kripke frames is easy to see using
singleton sets and successor sets. 

Using the two axioms $\mathtt{At}^n$ and $\mathtt{R}^n$, our main theorem is the following:
\begin{theorem}\label{thm:completeness-kripke-frames} Let $\Theta \subseteq
\mathcal{L}_{qf}$ be a set of Sahlqvist formulas s.t. the class $\KFr(\Theta)$ of Kripke
frames validating $\Theta$ has diversity $n$. Then $\mathsf{K}_\Pi \Theta
\mathtt{Bc}\mathtt{At}^n\mathtt{R}^n$ is sound and strongly complete for $\KFr(\Theta)$.
\end{theorem}
An outline of the proof is in order. First, we show that finite diversity means finite
depth, and the logic $\mathsf{K}_\Pi \Theta \mathtt{Bc}\mathtt{At}^n\mathtt{R}^n$
recognizes this. This means that $\Box^{\le n}$ simulates global modality well enough.
Second, to satisfy a consistent set of formulas, we take a point-generated general frame
$\mathbf{F}_a$ in the canonical saturated general frame where $a$ extends this consistent
set of formulas. Third, we keep only the `isolated' worlds in $\mathbf{F}_a$, i.e., those
worlds whose singleton is a proposition, and arrive at $\mathbf{F}_a^{at}$. This step does
not disturb the truth of formulas at the remaining worlds. Finally, we expand the
propositional domain of $\mathbf{F}_a^{at}$ to the full powerset. We have to show that
this step again keeps the truth value of all formulas in $\mathcal{L}$ unchanged, and also
the validity of $\Theta$. These two points only depend on that the propositional
domain of $\mathbf{F}_a^{at}$ contains all singletons, is closed under taking successor
sets, and contains all the duplicate classes of the underlying Kripke frame. Lemma
\ref{lem.property-domain} establishes these three properties.

Now we begin the proof. Fix a set $\Theta$ of Sahlqvist formulas with the diversity of $\KFr(\Theta)$
being $n$ and let $\Lambda = \mathsf{K}_\Pi \Theta \mathtt{Bc}\mathtt{At}^n\mathtt{R}^n$.
We first establish a logical point:
\begin{lemma}\label{lem:s4modality} $\Diamond^{n+1}p \to \Diamond^{\le n} p$ is a theorem
of $\Lambda$. Thus, denoting $\Diamond^{\le n}$ by $\sfe$ and the dual $\Box^{\le n}$ by
$\sfa$, $\Lambda$ proves that $\sfa$ is an $\mathtt{S4}$ modality that commutes with
$\forall$, and $\sfa$ works like the reflexive and transitive closure of $\Box$ in that
for example (1) $\sfe\Diamond\varphi \to \sfe\varphi \in \Lambda$ and (2) for any $m \in
\mathbb{N}$, $\Diamond^m\varphi \to \sfe\varphi \in \Lambda$.
\end{lemma}
\begin{proof}
First we show that any Kripke frame $\mathbf{G} \in \KFr(\Theta)$ must
also validate $\Diamond^{n+1}p \to \Diamond^{\le n} p$. Suppose not, then we have some $w
R x_1 R x_2 \dots R x_n R u$, which is also a shortest path from $w$ to $u$. This path is
also present and shortest in the subframe $\mathbf{G}_w$ of $\mathbf{G}$ generated from
$w$. Now note that $w, x_1, \dots, x_n$ are pairwise non-duplicates within $\mathbf{G}_w$,
since if there were a duplicate pair, then the path can be shortened. This contradicts
that $\KFr(\Theta)$ has diversity $n$. Hence any Kripke frame $\mathbf{G} \in \KFr(\Theta)$ validates $\Diamond^{n+1}p \to \Diamond^{\le n} p$. Since the normal propositional modal logic axiomatized by $\Theta$ is Kripke complete,
$\Diamond^{n+1} p \to \Diamond^{\le n} p$ is in $\Lambda$. The remaining claims follow
easily from basic normal modal reasoning and $\ttbc$.
\end{proof}
\noindent In the following we continue using $\sfa$ for $\Box^{\le n}$ and $\sfe$ for
$\Diamond^{\le n}$ and drop the superscripts on $\ttq^n$, $\ttat^n$, and $\ttr^n$. 

Now we start with the canonical saturated general frame $\bff_\Lambda = (W, R, B)$. Recall
that this involves expanding the language to $\mathcal{L}^+$ built from variables in
$\Prop^+$ and extending $\Lambda$ conservatively to $\Lambda^+$. Lemma
\ref{lem:s4modality} transfer to $\Lambda^+$ without problems. For Theorem
\ref{thm:completeness-kripke-frames}, clearly it is enough to show that every $w \in W$
can be satisfied in a Kripke frame validating $\Theta$. Thus fix an arbitrary $a \in W$
and consider the general frame $\bff_a = (W_a, R_a, B_a)$ generated from $a$, defined as
follows: 
\begin{itemize}
\item $W_a$ is $R^*[a]$ where $R^*$ is the reflexive and transitive closure of $R$;
\item $R_a = R \cap (W_a \times W_a)$;
\item $B_a = \{X \cap W_a \mid X \in B\}$; we write $[\varphi]_a$ for $[\varphi] \cap
W_a$, and with this notation, $B_a = \{[\varphi]_a \mid \varphi \in \mathcal{L}^+\}$.
\end{itemize} 

We show that $\sfa$ and $\sfe$ work as universal and existential modalities at $a$ and
$\ttq$ works as intended.
\begin{lemma}
  \label{lem:uni-exi-atom} 
  For any $\varphi \in \mathcal{L}^+$, $\sfa\varphi \in a$ iff $[\varphi]_a = W_a$, and
  similarly $\sfe\varphi \in a$ iff $[\varphi]_a$ is non-empty. Also, $\ttq(\varphi) \in a$
  iff $[\varphi]_a$ is a singleton.
\end{lemma}
\begin{proof}
For the first part, we just prove $\sfe\varphi \in a$ iff $[\varphi]_a \not= \varnothing$.
Since $a$ is an MCS and together with Lemma \ref{lem:s4modality}, $\sfe\varphi \in a$ iff
there is $m$ s.t. $\Diamond^m \varphi \in a$. By standard reasoning in canonical models,
i.e., repeated use of Lemma \ref{lem:existence-general}, this is true iff there is $u \in
W_a$ s.t. $\varphi \in u$. 

Now for $\ttq(\varphi)$, again since $a$ is a saturated MCS, $\ttq(\varphi) \in a$ iff (1)
$\sfe\varphi \in a$ and (2) for any $\psi \in \mathcal{L}^+$, $\sfa(\varphi \to \psi) \lor
\sfa(\varphi \to \lnot\psi) \in a$. (The second point uses saturation.) Using the first
part of this lemma, (1) translates to $[\varphi]_a \not= \varnothing$, and (2) translates
to that for any $\psi \in \mathcal{L}^+$, $[\varphi]_a \subseteq [\psi]_a$ or $[\varphi]_a
\subseteq (W_a \setminus [\psi]_a)$. If $[\varphi]_a$ is a singleton, these two points are
clearly true. Conversely, if $[\varphi]_a$ is empty, (1) is clearly false. If instead
there are distinct $x, y \in [\varphi]_a$, then there must be a formula $\psi \in
\mathcal{L}^+$ s.t. $\psi \in x$ but $\psi \not\in y$, making (2) false.
\end{proof}
This means the atoms of $B_a$ (as a Boolean algebra under set inclusion) are precisely
$\{[\varphi]_a \mid \mathtt{Q}(\varphi) \in a\}$. We focus on these atoms and define the
\emph{atomic subframe} $\bff_a^{at}$ of $\bff_a$ as $(W_a^{at}, R_a^{at}, B_a^{at})$ where 
\begin{itemize}
\item $W_a^{at} = \{w \in W_a \mid \{w\} \in B_a\}$;
\item $R^{at}_a = R_a \cap (W_a^{at} \times W_a^{at})$;
\item $B_a^{at} = \{X \cap W_a^{at} \mid X \in B_a\}$.
\end{itemize}
We write $[\varphi]_a^{at} = [\varphi]_a \cap W_a^{at}$. Then $B_a^{at} =
\{[\varphi]_a^{at} \mid \varphi \in \mathcal{L}^+\}$. A key property of $\bff_a^{at}$ is
that every world in it is named by a formula given how $W_a^{at}$ is defined. For each $w
\in W_a^{at}$, we fix a formula $\chi_w$ s.t. $\{w\} = [\chi_w]_a$. 

We want to immediately make sure that $a \in W_a^{at}$.
\begin{lemma}
The singleton $\{a\}$ is in $B_a$ and thus $a \in W_a^{at}$.
\end{lemma}
\begin{proof}
By a formal derivation, $\exists p (p \land \ttq(p))$ is in $\Lambda$. Indeed, with
$\mathtt{S4}$ normal modal reasoning and $\mathtt{Bc}$, we can derive $\ttq(p) \to
\mathsf{A}(p \to \ttq(p))$ in $\Lambda^+$. The main steps include 
  \begin{itemize}
\item $(\sfa(p \to q) \lor \sfa(p \to \lnot q)) \to \sfa(\sfa(p \to q) \lor \sfa(p \to
\lnot q))$
\item $\forall q(\sfa(p \to q) \lor \sfa(p \to \lnot q)) \to \sfa\forall q(\sfa(p \to q)
\lor \sfa(p \to \lnot q))$
\item $\sfa (p \to \sfe p)$ and then $\forall q(\sfa(p \to q) \lor \sfa(p \to \lnot q)) \to \sfa(p \to \ttq(p))$.
  \end{itemize}
Now suppose $\lnot \exists p (p \land \ttq(p))$. Then we derive $\sfe \forall p (p \to
\lnot \ttq(p))$. With $\ttat$, we derive $\exists p (\ttq(p) \land \sfa(p \to \forall p (p
\to \lnot\ttq(p))))$. A contradiction follows. Using (EE), we only need to derive a
contradiction from $\ttq(p) \land \sfa(p \to \forall p (p \to \lnot\ttq(p)))$. From
$\sfa(p \to \forall p(\ttq(p) \to \lnot p))$ we have $\sfa(p \to \lnot\ttq(p))$. And
recall $\ttq(p) \to \sfa(p \to \ttq(p))$ is in $\Lambda^+$. This means we derive $\sfa(p
\to \ttq(p))$ and then $\sfa \lnot p$, which contradicts the $\sfe p$ part in $\ttq(p)$.

Thus $\exists p (p \land \ttq(p)) \in a$, and since $a$ is saturated, there is $r \in
\Prop^+$ s.t. $r \land \ttq(r) \in a$. Then $[r]_a$ must be $\{a\}$ and $\{a\} \in B_a$.
\end{proof}

From this, we show that $\bff_a^{at}$ behaves as a canonical general frame and thus is
quantifiable:
\begin{lemma}\label{lem:existence-truth-atomic-frame}
  \begin{itemize}
\item For any $w \in W_a^{at}$ and $\varphi \in \mathcal{L}^+$, $\varphi \in w$ iff
$\sfa(\chi_w \to \varphi) \in a$ iff $\sfe(\chi_w \land \varphi) \in a$. Also, for any $w,
u \in W_a^{at}$, $wR_a^{at}u$ iff $\sfe(\chi_w \land \Diamond \chi_u) \in a$.
\item For any $\Diamond\varphi \in w \in W_a^{at}$, there is $u \in R_a^{at}[w]$ s.t.
$\varphi \in u$.
\item For any substitution $\sigma$, define the associated valuation $[\sigma]_a^{at}$ for
$\bff_a^{at}$ by $[\sigma]_a^{at}(p) = [\sigma(p)]_a^{at}$. Then
$\tset{\varphi}^{\bff_a^{at}}([\sigma]_a^{at}) = [\sigma(\varphi)]_a^{at}$.
  \end{itemize}
\end{lemma}
\begin{proof}
By definition, $[\chi_w]_a = \{w\}$ and $[\chi_u]_a = \{u\}$. As we have shown,
$\sfe(\chi_w \land \Diamond\chi_u) \in a$ iff $[\chi_w \land \Diamond\chi_u]_a$ is
non-empty, iff $\Diamond\chi_u \in w$, and iff $w R_a^{at} u$. The universal case is dual to the existential case.
  
Now suppose $\Diamond\varphi \in w \in W_a^{at}$. First, apply $\ttr$ to $\chi_w$.
Using saturation at $a$, we have a formula (indeed a variable) which we denote by
$\chi_{R[w]}$ s.t. $\sfa(\chi_w \to \Box\chi_{R[w]}) \land \forall r (\sfa(\chi_w \to \Box
r) \to \sfa(\chi_{R[w]} \to r))$ is in $a$. By plugging in $\lnot r$ for $r$ and
contraposing, we have $\forall r (\sfe(\chi_{R[w]} \land r) \to \sfe(\chi_w \land \Diamond
r)) \in a$. Also, from $\sfa(\chi_w \to \Box\chi_{R[w]}) \in a$ and $\Diamond\varphi \in
w$, $\sfe(\chi_w \land \Diamond(\chi_{R[w]} \land \varphi)) \in a$. This means
$\sfe\Diamond(\chi_{R[w]} \land \varphi)$ and thus $\sfe(\chi_{R[w]} \land \varphi) \in
a$. Then by $\ttat$ (or the atomicity of $B_a$), there is $u \in W_a^{at}$ s.t. $\varphi$
and $\chi_{R[w]} \in u$. Then $\sfe(\chi_{R[w]} \land \chi_u) \in a$, and thus
$\sfe(\chi_w \land \Diamond \chi_u) \in a$, which means $u \in R_a^{at}[w]$. 

For the last part, recall that for the original canonical saturated frame
$\bff_{\Lambda}$, Lemma \ref{lem:truth-substitution} applies: for any $\varphi$ and
substitution $\sigma$ for $\mathcal{L}^+$, $\tset{\varphi}^{\bff_\Lambda}([\sigma]) =
[\sigma(\varphi)]$. $\bff_a^{at}$ is obtained by restricting $\bff_{\Lambda}$ to
$W_{a}^{at}$. Thus, it is enough to show that for any valuation $v$ for $\bff_\Lambda$,
writing $v|_a^{at}$ for the restricted valuation defined by $v|_a^{at}(p) = v(p) \cap
W_a^{at}$, $\tset{\varphi}^{\bff_a^{at}}(v|_a^{at}) = \tset{\varphi}^{\bff_\Lambda}(v)
\cap W_a^{at}$. Using induction on $\varphi$, the base and the Boolean cases are trivial,
as the operation of relative negation and intersection commutes intersecting with
$W_a^{at}$. For the modal case, we need $m_{\Diamond}^{\bff_a^{at}}(X \cap W_a^{at}) =
m_{\Diamond}^{\bff_\Lambda}(X) \cap W_a^{at}$ where $X$ is assumed to be in $B$ due to IH.
The left-to-right inclusion is trivial since $\bff_a^{at}$ is a restriction of
$\bff_{\Lambda}$. For the right-to-left inclusion, first write $X$ as $[\psi]$ for some
$\psi \in \mathcal{L}^+$ and use the second bullet point of this lemma. The quantifier
case is not much different from the case for disjunction, using only the distribution of
intersection over arbitrary union and that $B_a^{at} = \{X \cap W_a^{at} \mid X \in B\}$.
\end{proof}
Now we start to show that $B_a^{at}$ has what it takes to be a $\pqml$-invariant subdomain of
$\wp(W_a^{at})$ over the underlying Kripke frame $\bbf_a^{at} = (W_a^{at}, R_a^{at})$.
\begin{lemma}\label{lem.property-domain}
  \begin{itemize}
\item For any $w \in W_a^{at}$, $\{w\} \in B_a^{at}$.
\item For any $X \in B_a^{at}$, $R_a^{at}[X] \in B_a^{at}$.
\item $\bff_a^{at}$ is point-generated from $a$ and $\bbf_a^{at} \models \Theta$. Thus it
has diversity $k \le n$.
\item Let $\Delta$ be the duplicate relation of $\bbf_a^{at}$. Then each $D \in
W_a^{at}/\Delta$ is in $B_a^{at}$.
  \end{itemize}
\end{lemma}
\begin{proof}
The first bullet point is trivial. For the second bullet point, pick any $X \in B_a^{at}$.
Then we have some $\varphi$ s.t. $X = [\varphi]_a^{at}$. Reasoning in the saturated MCS
$a$ and apply $\ttr$ to $\varphi$, we obtain a $\psi$ s.t. $\sfa(\varphi \to \Box\psi) \in
a$ and for any $\gamma \in \mathcal{L}^+$, $\sfe(\psi \land \gamma) \to \sfe(\varphi \land
\Diamond\gamma) \in a$. Now for any $w \in X$, $\varphi \in w \in W_a^{at}$. Thus
$\sfe(\chi_w \land \varphi) \in a$. Together with $\sfa(\varphi \to \Box\psi)$,
$\sfe(\chi_w \land \Box\psi) \in a$, meaning $\Box\psi \in w$. Thus $R[w] \subseteq
[\psi]$ and hence $R_a^{at}[w] \subseteq [\psi]_a^{at}$. Since $w$ is chosen arbitrarily
from $X$, $R_a^{at}[X]\subseteq [\psi]_a^{at}$. On the other hand, suppose $u \in
[\psi]_a^{at}$. Then $\sfe(\psi \land \chi_u) \in a$. Then $\sfe(\varphi \land
\Diamond\chi_u) \in a$ and thus there is $w \in W_a^{at}$ s.t. $\varphi$ and
$\Diamond\chi_u$ are in $w$. This means $w \in X$ and $w R_a^{at} u$. So in sum,
$[\psi]_a^{at} \subseteq R_a^{at}[X]$.

For the third bullet point, we need to first show that every $w \in W_a^{at}$ is reachable
from $a$ within $\bff_a^{at}$. If $w \in W_a^{at}$, then at least $\sfe\chi_w \in a$. This
means for some $m \le n$, $\Diamond^m\chi_w \in a$. By repeated use of the second bullet
point of Lemma \ref{lem:existence-truth-atomic-frame}, there is indeed a path from $a$ to
$w$ inside $\bff_a^{at}$. Now, by the third bullet point of Lemma
\ref{lem:existence-truth-atomic-frame}, $\bff_a^{at} \models \Theta$. For this to transfer
to $\bbf_a^{at}$, we rely on the assumption that $\Theta$ consists of Sahlqvist formulas,
which are $\mathcal{AT}$-persistent in the sense that if they are valid on a general frame
whose set of admissible sets contains all singletons (the atomic/$\mathcal{A}$ part) and
is closed under taking successor set (the tense/$\mathcal{T}$ part), then they are also
valid on the underlying Kripke frame. That Sahlqvist formulas are
$\mathcal{AT}$-persistent has been observed for example in \cite{venema1993derivation}.
The idea is that for any $\varphi \in \Theta$, if $\bbf_a^{at} \not\models \varphi$, then
there is a falsifying valuation that only uses sets in $B_a^{at}$ so that it is also a
valuation for $\bff_a^{at}$, contradicting that $\bff_a^{at} \models \Theta$. This special
valuation is obtained by the standard minimal valuation technique for Sahlqvist formulas
and note that in minimal valuations, only finite unions of sets of the form
$(R_a^{at})^m[\{w\}]$ are used, which are in $B_{a}^{at}$ by the assumed closure
properties.

Now that $\bbf_a^{at}$ has diversity $k \le n$, there are also $b_1, b_2, \dots, b_k \in
W_a^{at}$ each representing a duplicate class. For each $i$, we show that the duplicate
class $D$ that $b_i$ is in is in $B_a^{at}$. Now if $D = \{b_i\}$ then we have shown that
it is in $B_a^{at}$. So assume that there is $c \not= b_i$ in $D$. Observe that for any $w
\in W_a^{at} \setminus \{b_1, \dots, b_k\}$, $w \in D$ iff $w$ and $c$ are duplicates, and
iff the following are true:
  \begin{itemize}
\item for any $j$, $w R_a^{at} b_j$ iff $c R_a^{at} b_j$;
\item for any $j$, $b_j R_a^{at} w$ iff $b_j R_a^{at} c$;
\item $w$ is reflexive iff $c$ is reflexive;
\item $w R_a^{at} c$ iff $c R_a^{at} w$.
  \end{itemize}
The above conditions are all expressible in $B_a^{at}$ using singletons, the
$R_a^{at}[\cdot]$ operation, the $m_{\Diamond}^{\bbf_a^{at}}$ operation, and also the set
of reflexive points in $\bbf_a^{at}$ defined by sentence $\forall p (\Box p \to p)$. In
fact, in the original canonical saturated general frame $\bbf_{\Lambda}$, $[\forall p
(\Box p \to p)]$ is already the set of all reflexive worlds since any two worlds are
separated by a proposition in $B$. Indeed, the set
  \begin{align*}
    \{m_\Diamond^{\bbf_a^{at}}(\{b_j\}) \mid c R_a^{at} b_j\} \cup 
    \{W_a^{at} \setminus m_\Diamond^{\bbf_a^{at}}(\{b_j\}) \mid \text{not } c R_a^{at} b_j\} \cup  \\
    \{R_a^{at}[b_j] \mid b_j R_a^{at} c\} \cup 
    \{W_a^{at} \setminus R_a^{at}(b_j) \mid \text{not } b_j R_a^{at} c\} \cup \\
    \{[\forall p (\Box p \to p)]_a^{at} \mid c R_a^{at} c\} \cup 
    \{W_a^{at} \setminus [\forall p (\Box p \to p)]_a^{at} \mid \text{not } c R_a^{at} c\} \cup \\
    \{(m_\Diamond^{\bbf_a^{at}}(\{c\}) \cap R_a^{at}[c]) \cup ((W_a^{at} \setminus m_\Diamond^{\bbf_a^{at}}(\{c\})) \cap (W_a^{at} \setminus R_a^{at}[c]))\}
  \end{align*}
contains all the required conditions, the intersection of which we denote by $X$. Then $(X \setminus \{b_1, \dots, b_k\}) \cup \{b_i\}$ is the duplicate class $D$ that $b_i$ is in.
\end{proof}

Putting pieces together, for any $\Lambda^+$-MCS $\Sigma$, it is satisfied on
$\bff^{at}_a$ by itself under valuation $[\iota]^{at}_a$ by Lemma
\ref{lem:existence-truth-atomic-frame}. By Lemma \ref{lem.property-domain} and Theorem
\ref{thm:ele-subdomain}, $B^{at}_a$ is a $\pqml$-invariant subdomain of $\wp(W^{at}_a)$. By
definition of $\pqml$-invariant subdomain and $[\iota]^{at}_a\in (B^{at}_a )^{\Prop}$,
$\Sigma$ is satisfied on $\bbf^{at}_a$ by itself under valuation $[\iota]^{at}_a$.
Finally, by the third bullet point of Lemma \ref{lem.property-domain}, $\bbf^{at}_a$ is a
frame validating $\Theta$. This completes the proof of theorem
\ref{thm:completeness-kripke-frames}.

By observing where we used the Sahlqvist condition on $\Theta$, we note a different way
of stating our main result.

\begin{definition}
For any formula $\varphi \in \mathcal{L}$, it is \emph{$\mathcal{ATQ}$-persistent} if for
any quantifiable general frame $(W, R, B)$ such that all singleton subsets of $W$ is in
$B$ and $B$ is closed under $R[\cdot]$, if $\varphi$ is valid on $(W, R, B)$, then
$\varphi$ is valid on the underlying $(W, R)$.
\end{definition}
\begin{definition}
For any \npqml\ $\Lambda \subseteq \mathcal{L}$, it is Kripke $\mathcal{L}_{qf}$-complete
if for any $\varphi \in \mathcal{L}_{qf}$ that is valid on all Kripke frames in
$\mathsf{KFr}(\Lambda)$, $\varphi$ is in $\Lambda$.
\end{definition}
\begin{corollary}
Let $\Theta \subseteq \mathcal{L}$ be a set of $\mathcal{ATQ}$-persistent formulas such
that $\mathsf{KFr}(\Theta)$ has diversity $n$ and $\mathsf{K}_\Pi \Theta
\mathtt{Bc}\mathtt{At}^n\mathtt{R}^n$ is Kripke $\mathcal{L}_{qf}$-complete. Then
$\mathsf{K}_\Pi \Theta \mathtt{Bc}\mathtt{At}^n\mathtt{R}^n$ is the \npqml\ of
$\mathsf{KFr}(\Theta)$.
\end{corollary}

Now we consider a special case, the class $\KFr(\ttt{5})$ of Euclidean frames. We have
observed that $\KFr(\mathtt{5})$ has diversity $3$. Note that it also validates
$\Diamond^3 p \to \Diamond^2 p$. Thus, the \npqml\ of Euclidean Kripke frames can be
axiomatized as $\mathsf{K}_\Pi \mathtt{5} \mathtt{Bc}\mathtt{At}^2\mathtt{R}^2$. In other
words, $\ttt{5}\pi+$ is $\mathsf{K}_\Pi \mathtt{5} \mathtt{Bc}\mathtt{At}^2\mathtt{R}^2$.
However, unlike $\ttt{D45}\pi+ = \mathsf{K}_\Pi\mathtt{D45}\mathtt{Bc}\mathtt{At}^1$,
$\mathtt{R}^2$ is indispensable for the axiomatization of the logic of Euclidean frames.
For this, we construct a quantifiable frame $\mathbf{F}$ such that $\mathbf{F}$ validates
$\mathsf{K}_\Pi\ttt{5BcAt}^2$ but not $\exists p (\Box p \land \forall q (\Box q \to
\Box^2(p \to q)))$, a formula that is valid on any Kripke frame using the successor sets.

Let $\bff = (W,R,B)$ where
\begin{itemize}
  \item $W = \mathbb{N}\cup \{w\}$, where $w \not \in \mathbb{N}$;
  \item $R = \mathbb{N}^2\cup \{(w,2n)\mid n\in \mathbb{N}\}$; i.e. $R$ is total in $\mathbb{N}$ and $w$ sees all even numbers;
  \item $B = \{X\subseteq W\mid X \text{ is finite or cofinite with respect to $W$}\}$.
\end{itemize}

\begin{proposition}
    $\bff$ is a quantifiable frame validating $\mathsf{K}_\Pi \mathtt{5BcAt^2}$ but not 
    $\exists p(\Box p\land \forall q(\Box q\to \Box^2(p\to q)))$.
\end{proposition}
\begin{proof}
First, let us check that $\bff$ is a quantifiable frame. Consider first the subframes
$\bfg = (\NN, \NN^2, B')$ generated from (any) $n \in \NN$, where $B' = \{X \cap \NN \mid X
\in B\}$. Observe that $B'$ is also $\{X \subseteq \NN \mid X \text{ is finite or } \NN
\setminus X \text{ is finite}\}$. We first show that $\bfg$ is
quantifiable. There is only one duplicate class $\NN$ of $(\NN, \NN^2)$, and $B'$ is a
discrete field of sets, so Theorem \ref{thm:ele-subdomain} applies: for any formula
$\varphi$ and valuation $v: \mathsf{Prop} \to B'$, $\tset{\varphi}^{B'}(v) =
\tset{\varphi}^{\wp(\NN)}(v)$. By Lemma \ref{lem:boolean-breakdown}, there is then a
Boolean formula $\alpha$ such that $\tset{\varphi}^{\wp(\NN)}(v) = \tset{\alpha}(v)$. So
$\tset{\varphi}^{B'}(v)$ is a Boolean combination of propositions in $B'$, which is in
$B'$ as $B'$ is a field of sets. This shows that $\bfg$ is quantifiable.

Now, since for any $n \in \NN$, its point-generated subframe is $\bfg$, by an easy
induction, for any valuation $v$ (with $v|_\NN$ being its point-wise restriction to $\NN$)
and any formula $\varphi$, $n \in \tset{\varphi}^\bff(v)$ iff $n \in
\tset{\varphi}^\bfg(v|_\NN)$. In other words, $\tset{\varphi}^\bff(v) \cap \NN =
\tset{\varphi}^\bfg(v|_\NN)$. As $\bfg$ is quantifiable, $\tset{\varphi}^\bff(v) \cap \NN$
is finite or cofinite in $\NN$. But $\tset{\varphi}^\bff(v)$ and $\tset{\varphi}^\bff(v)
\cap \NN$ differ by at most one point $w$. So $\tset{\varphi}^\bff(v)$ is finite or
cofinite in $W$, and is in $B$. This shows that $\bff$ is quantifiable.

As $\bff$ is quantifiable, as soon as we verify that $\bff$ validates $\ttt{5}$ and
$\ttt{At}^3$, $\bff$ then automatically validates $\mathsf{K}_\Pi\ttt{5BcAt}^3$. But
the underlying Kripke frame of $\bff$ already validates $\ttt{5}$, and $B$ contains
all singleton subsets of $W$. 

To see that $\exists p(\Box p\land \forall q(\Box q\to \Box^2(p\to q)))$ is invalid, we
only need to note that for any $X \in B$ such that $X \supseteq R[w]$, there is always a
$Y \in B$ such that $Y \subsetneq X$ and $Y \supseteq R[w]$. This is because $X$ must
contain infinitely many odd numbers, and we can delete one odd number to obtain $Y$. This
shows that at $w$, the formula $\forall p (\Box p \to \exists q (\Box q \land \Diamond^2(p
\land \lnot q)))$ is true at $w$. Thus, its negation $\exists p (\Box p \land \forall q
(\Box q \to \Box^2(p \to q)))$ is false at $w$.
\end{proof}
Thus, $\mathsf{K}_\Pi\ttt{5Bc}\ttt{At}^2$ is not the logic of Euclidean frames.

We can also show that the requirement of $\Theta$ consisting of Sahlqvist formulas cannot
be dispensed with altogether either. Consider the axioms $\mathtt{T}$: $p\to \Diamond p$,
$\mathtt{M}$: $\Diamond \Box \neg p\vee \Diamond \Box p$, $\mathtt{E}$: $\Diamond
(\Diamond p\land \Box q)\to \Box (\Diamond p\vee \Box q)$, and $\mathtt{Q}$: $(\Diamond
p\land \Box (p\to \Box p))\to p$ (we reuse the letter \ttt{Q}) used in \cite{van1978two}.
It is not hard to show that the only Kripke frames validating \ttt{TMEQ} are those with
the identity accessibility relation, and thus $\KFr(\ttt{TMEQ})$ has diversity $1$.
However, $\mathsf{K}_\Pi \mathtt{TMEQ} \mathtt{Bc}\mathtt{At}^1\mathtt{R}^1$ is not the
\npqml\ of $\KFr(\ttt{TMEQ})$, as it does not derive $p \leftrightarrow \Box p$ that is
valid on Kripke frames with the identity accessibility relation. The idea is to observe
that the veild recession general frame $\bff = (\mathbb{Z}, R, B)$ where 
\begin{itemize}
  \item $nRm$ iff $m \ge n-1$, and 
\item $X \in B$ iff there is $n$ such that for all $m \ge n$, $m \in X$ iff $n \in X$
(call this property \emph{settled after $n$}),
\end{itemize}
is quantifiable and validates $\ttt{TMEQAt}^1\ttt{R}^1$, but not $p \to \Box p$. Thus, the
set of validities of $\bff$ separates $\mathsf{K}_\Pi\ttt{TMEQAt}^1\ttt{R}^1$ from $p \to
\Box p$ and shows that $\mathsf{K}_\Pi\ttt{TMEQAt}^1\ttt{R}^1$ is not $\ttt{TMEQ}\pi+$. It
is easy to check that $\bff$ validates the special axioms, but it takes more work to show
that $\bff$ is quantifiable. Here we need to use truncated point-generated submodels and
the key step is to show that for all sufficiently large $m$ and $n$, if $\exists p
\varphi$ is true at $m$ using set $X$ as the valuation for $p$ and $d$ is the modal depth
of $\exists p \varphi$, then $\exists p \varphi$ is also true at $n$ as we can shift $X$
by $n-m$ and then the depth-$d$ truncated submodel generated from $n$ is isomorphic to the
depth-$d$ truncated submodel generated from $m$, making this shifted $X$ a witness to
$\exists p \varphi$ for $n$. See appendix at the end of the paper for detailed proofs.

\section{Conclusion}\label{sec:conclusion}

We conclude with some ideas for possible future research. First, the result that all
\npqml s containing \ttt{Bc} is complete for the class of quantifiable frames it defines
is fairly standard and expected, but can still be generalized, for example, to
\emph{neighborhood} frames. Of course, \ttt{Bc} will be dropped from the logic, and
similar work has been done for first-order modal logic \cite{arlo2006first}. One may also
consider the alternative semantics in \cite{goldblatt2006general} where $\forall p
\varphi$ is true at $w$ iff there is a proposition $X$ that contains $w$ and entails all
propositions expressible by $\varphi$ as we vary the proposition denoted by $p$. We also
believe that it would be instructive to rewrite the proof of Theorem
\ref{thm:completeness-quantifiable-frames} in terms of Lindenbaum algebras and duality
theory, as this may help us generalize the result.

For our second result, there are two natural ideas to generalize. The first is dropping
$\mathtt{At}^n$ and consider completeness w.r.t. algebraic semantics based on complete and
completely multiplicative modal algebras. It should be noted that, as is clear in our
proof, axiom $\mathtt{At}^n$ corresponds to the existence of `world propositions' that can
later be interpreted as possible worlds, and the world propositions serve as the names of
the possible worlds. Hybrid logics use world propositions in a much more direct way by
taking them as a primitive syntactical category, namely the nominals, and a recent work
\cite{blackburn2023axiom} has considered propositionally quantified hybrid modal logic. As
is mentioned in that paper, Arthur Prior is a strong proponent of both. However, the idea
of there being no maximally specified possible worlds but only partial states
\cite{humberstone1981worlds,rumfitt2015boundary,Holliday2021-HOLPS-4,Ding2020-DINAPI} is
also worth investigating in this context (though \pqml\ together with plural quantifiers
are used to argue for there being world propositions \cite{fritz2023foundations}), and
algebraic semantics allowing atomless elements in the algebras is a natural way to model
this. Previous works in this line include \cite{Holliday2019note,Ding2018-DINOTL,Ding2021}.

The other direction for generalization is dropping the finite diversity condition in some
way. The condition that $\KFr(\Theta)$ has finite diversity is admittedly a very
restrictive one, and it is worth investigating the exact scope of this condition,
especially together with the requirement of $\Theta$ consisting of only Sahlqvist
formulas. We see that there is at least one promising way of relaxing the finite diversity
condition: requiring only finite diversity for each point-generated frame of finite depth. 

Finally, we mention a broader question: can the theory of \pqml s inform the theory of
modal $\mu$-calculus \cite{bradfield2007modal} or vice versa, especially over completeness
questions, noting that the $\mu$ operator is also a kind of propositional quantifier? For
example, the recent work \cite{FritzConservativity} utilized the fixpoint construction in
\npqml , and we believe more needs to be done.

\section*{Acknowledgements}
We thank the anonymous reviewers for their helpful suggestions that led to many
improvements. We also thank Peter Fritz for commenting on a very early draft of this
paper. The first author is supported by NSSF grant 22CZX066.


\newcommand{\van}[1]{}
\bibliographystyle{aiml}
\bibliography{aiml18}

\begin{thebibliography}{10}
\expandafter\ifx\csname url\endcsname\relax
  \def\url#1{\texttt{#1}}\fi
\expandafter\ifx\csname urlprefix\endcsname\relax\def\urlprefix{URL }\fi
\newcommand{\enquote}[1]{``#1''}

\bibitem{antonelli2002representability}
Antonelli, G.~A. and R.~H. Thomason, \emph{Representability in second-order
  propositional poly-modal logic}, The Journal of Symbolic Logic \textbf{67}
  (2002), pp.~1039--1054.

\bibitem{arlo2006first}
Arl{\'o}-Costa, H. and E.~Pacuit, \emph{First-order classical modal logic},
  Studia Logica \textbf{84} (2006), pp.~171--210.

\bibitem{bednarczyk2022whydoes}
Bednarczyk, B. and S.~Demri, \emph{Why does propositional quantification make
  modal and temporal logics on trees robustly hard?}, Logical Methods in
  Computer Science \textbf{18} (2022).

\bibitem{belardinelli2018secondAIJ}
Belardinelli, F., W.~Van Der~Hoek and L.~B. Kuijer, \emph{Second-order
  propositional modal logic: Expressiveness and completeness results},
  Artificial Intelligence \textbf{263} (2018), pp.~3--45.

\bibitem{bilkova2007uniform}
B{\'\i}lkov{\'a}, M., \emph{Uniform interpolation and propositional quantifiers
  in modal logics}, Studia Logica \textbf{85} (2007), pp.~1--31.

\bibitem{blackburn2023axiom}
Blackburn, P., T.~Bra{\"u}ner and J.~L. Kofod, \emph{An axiom system for basic
  hybrid logic with propositional quantifiers}, in: \emph{International
  Workshop on Logic, Language, Information, and Computation}, Springer, 2023,
  pp. 118--134.

\bibitem{bradfield2007modal}
Bradfield, J. and C.~Stirling, \emph{Modal mu-calculi}, in: \emph{Handbook of
  Modal Logic}, Elsevier, 2007 pp. 721--756.

\bibitem{ten2006expressivity}
{\van{Cate}}ten~Cate, B., \emph{Expressivity of second order propositional
  modal logic}, Journal of Philosophical Logic \textbf{35} (2006),
  pp.~209--223.

\bibitem{Dagostino2000logical}
D'agostino, G. and M.~Hollenberg, \emph{Logical questions concerning the
  $\mu$-calculus: Interpolation, {Lyndon} and {{\L}o{\'s}-Tarski}}, The Journal
  of Symbolic Logic \textbf{65} (2000), pp.~310--332.

\bibitem{Dekker2024KD45}
Dekker, P.~M., \emph{{KD45} with propositional quantifiers}, Logic and Logical
  Philosophy \textbf{33} (2024), pp.~27--54.

\bibitem{Ding2018-DINOTL}
Ding, Y., \emph{On the logics with propositional quantifiers extending
  {S5}\ensuremath{\Pi}}, in: G.~Bezhanishvili, G.~D'Agostino, G.~Metcalfe and
  T.~Studer, editors, \emph{Advances in Modal Logic, Vol. 12}, College
  Publications, 2018 pp. 219--235.

\bibitem{Ding2021}
Ding, Y., \emph{On the logic of belief and propositional quantification},
  Journal of Philosophical Logic \textbf{50} (2021), pp.~1143--1198.

\bibitem{Ding2020-DINAPI}
Ding, Y. and W.~H. Holliday, \emph{Another problem in possible world
  semantics}, in: N.~Olivetti and R.~Verbrugge, editors, \emph{Advances in
  Modal Logic, Vol. 13}, College Publications, 2020 pp. 149--168.

\bibitem{Dagostino2005axiomatization}
D’Agostino, G. and G.~Lenzi, \emph{An axiomatization of bisimulation
  quantifiers via the $\mu$-calculus}, Theoretical Computer Science
  \textbf{338} (2005), pp.~64--95.

\bibitem{fine1970propositional}
Fine, K., \emph{Propositional quantifiers in modal logic}, Theoria \textbf{36}
  (1970), pp.~336--346.

\bibitem{French2006thesis}
French, T., \enquote{Bisimulation quantifiers for modal logics,} {PhD} thesis,
  The University of Western Australia (2006).

\bibitem{french2007idempotent}
French, T., \emph{Idempotent transductions for modal logics}, in:
  \emph{Frontiers of Combining Systems: 6th International Symposium, FroCoS
  2007 Liverpool, UK, September 10-12, 2007 Proceedings 6}, Springer, 2007, pp.
  178--192.

\bibitem{French2002AiML}
French, T. and M.~Reynolds, \emph{A sound and complete proof system for
  {QPTL}}, in: P.~Balbiani, N.-Y. Suzuki, F.~Wolter and M.~Zakharyaschev,
  editors, \emph{Advances in Modal Logic, Vol. 4}, King's College Publications,
  2002 pp. 127--148.

\bibitem{fritz2022}
Fritz, P., \emph{Axiomatizability of propositionally quantified modal logics on
  relational frames}, The Journal of Symbolic Logic  (2022), p.~1–36.

\bibitem{fritz2023foundations}
Fritz, P., \enquote{The Foundations of Modality: From Propositions to Possible
  Worlds,} Oxford University Press, 2023.

\bibitem{FritzConservativity}
Fritz, P., \emph{Nonconservative extensions by propositional quantifiers and
  modal incompleteness}, manuscript  (2023).

\bibitem{Fritz2024}
Fritz, P., \enquote{Propositional Quantifiers,} Elements in Philosophy and
  Logic, Cambridge University Press, 2024.

\bibitem{ghilardi1995undefinability}
Ghilardi, S. and M.~Zawadowski, \emph{Undefinability of propositional
  quantifiers in the modal system {S4}}, Studia Logica \textbf{55} (1995),
  pp.~259--271.

\bibitem{goldblatt2006general}
Goldblatt, R. and E.~D. Mares, \emph{A general semantics for quantified modal
  logic.}, in: G.~Governatori, I.~Hodkinson and Y.~Venema, editors,
  \emph{Advances in modal logic, Vol. 6}, College Publications, 2006 pp.
  227--246.

\bibitem{Holliday2019note}
Holliday, W.~H., \emph{A note on algebraic semantics for {S5} with
  propositional quantifiers}, Notre Dame Journal of Formal Logic \textbf{60}
  (2019), pp.~311--332.

\bibitem{Holliday2021-HOLPS-4}
Holliday, W.~H., \emph{Possibility semantics}, in: M.~Fitting, editor,
  \emph{Selected Topics from Contemporary Logics}, College Publications, 2021
  pp. 363--476.

\bibitem{humberstone1981worlds}
Humberstone, I.~L., \emph{From worlds to possibilities}, Journal of
  Philosophical Logic  (1981), pp.~313--339.

\bibitem{kaminski1996expressive}
Kaminski, M. and M.~Tiomkin, \emph{The expressive power of second-order
  propositional modal logic}, Notre Dame Journal of Formal Logic \textbf{37}
  (1996), pp.~35--43.

\bibitem{kremer1997complexity}
Kremer, P., \emph{On the complexity of propositional quantification in
  intuitionistic logic}, The Journal of Symbolic Logic \textbf{62} (1997),
  pp.~529--544.

\bibitem{kremer2018completeness}
Kremer, P., \emph{Completeness of second-order propositional {S4} and {H} in
  topological semantics}, The Review of Symbolic Logic \textbf{11} (2018),
  pp.~507--518.

\bibitem{kuhn2004simple}
Kuhn, S. et~al., \emph{A simple embedding of {T} into double {S5}}, Notre Dame
  Journal of Formal Logic \textbf{45} (2004), pp.~13--18.

\bibitem{kuusisto2015second}
Kuusisto, A., \emph{Second-order propositional modal logic and monadic
  alternation hierarchies}, Annals of Pure and Applied Logic \textbf{166}
  (2015), pp.~1--28.

\bibitem{Cresswell1996-CREANI-3}
M.~J.~Cresswell, G. E.~H., \enquote{A New Introduction to Modal Logic,}
  Routledge, New York, 1996.

\bibitem{pitts1992interpretation}
Pitts, A.~M., \emph{On an interpretation of second order quantification in
  first order intuitionistic propositional logic}, The Journal of Symbolic
  Logic \textbf{57} (1992), pp.~33--52.

\bibitem{riba2012model}
Riba, C., \emph{A model theoretic proof of completeness of an axiomatization of
  monadic second-order logic on infinite words}, in: \emph{Theoretical Computer
  Science: 7th IFIP TC 1/WG 2.2 International Conference, TCS 2012, Amsterdam,
  The Netherlands, September 26-28, 2012. Proceedings 7}, Springer, 2012, pp.
  310--324.

\bibitem{rumfitt2015boundary}
Rumfitt, I., \enquote{The boundary stones of thought: An essay in the
  philosophy of logic,} Oxford University Press, USA, 2015.

\bibitem{steinsvold2020some}
Steinsvold, C., \emph{Some formal semantics for epistemic modesty}, Logic and
  Logical Philosophy \textbf{29} (2020), pp.~381--413.

\bibitem{van1978two}
\van{Benthem}van Benthem, J., \emph{Two simple incomplete modal logics},
  Theoria \textbf{44} (1978), pp.~25--37.

\bibitem{ten2004model}
\van{Cate}ten Cate, B., \enquote{Model theory for extended modal languages,}
  {PhD} thesis, University of Amsterdam (2004).

\bibitem{venema1993derivation}
Venema, Y., \emph{Derivation rules as anti-axioms in modal logic}, The Journal
  of Symbolic Logic \textbf{58} (1993), pp.~1003--1034.

\bibitem{visser1996bisimulations}
Visser, A., \emph{Bisimulations, model descriptions and propositinal
  quantifiers}, Logic Group Preprint Series \textbf{161} (1996).

\bibitem{zach2004decidability}
Zach, R., \emph{Decidability of quantified propositional intuitionistic logic
  and {S4} on trees of height and arity $\le\omega$}, Journal of Philosophical
  Logic \textbf{33} (2004), pp.~155--164.

\end{thebibliography}

\section*{Appendix}

We show that in our result Theorem \ref{thm:completeness-kripke-frames} the Sahlqvist
condition cannot be dropped completely. Without it, Kripke incompleteness for the
quantifier-free fragment can lead to incompleteness for the full propositionally
quantified language even with $\mathtt{At}^n$ and $\mathtt{R}^n$. Recall the axioms:
\begin{itemize}
\item $\mathtt{T}: \Box p\to p$;
\item $\mathtt{M}: \Box \Diamond p\to \Diamond \Box p$;
\item $\mathtt{E}: \Diamond(\Diamond p\land \Box q)\to \Box (\Diamond p\vee \Box q)$;
\item $\mathtt{Q}: (\Diamond p\land \Box (p\to \Box p))\to p$. 
\end{itemize}
We know that $\mathtt{TMEQ}$ defines the class of Kripke frames with the identity
accessibility relation and is Kripke incomplete \cite{van1978two}. Thus,
the diversity of $\KFr(\ttt{TMEQ})$ is just $1$. We will show the following:
\begin{theorem}
$\mathsf{K}_\Pi\mathtt{TMEQ} \mathtt{Bc}\mathtt{At}^1\mathtt{R}^1\neq
\mathtt{TMEQ}\pi+$.
\end{theorem}
The idea is to find a quantifiable general frame that validates \ttt{TMEQ} and $\ttt{At}^1$ and $\ttt{R}^1$, but not $p \to \Box p$, which is valid in $\KFr(\ttt{TMEQ})$. To facilitate the proofs, we first define a more liberal version of general frames (not requiring any closure conditions for the propositional domain) and then prove a lemma on how the truth of formulas is preserved under taking truncated generated submodels.



\begin{definition}
A \emph{frame with a propositional domain} (\emph{pd-frame} for short) is a triple $(W, R, B)$ where $(W, R)$ is a Kripke frame and $B \subseteq \wp(W)$ is non-empty. (This is only dropping the closure properties for $B$ in general frames.) We often conflate $\bff = (W, R, B)$ with $W$ when no confusion would arise. Valuations and semantics for pd-frames are defined exactly the same as for general frames. 

A \emph{pd-model} is a pair $\mathcal{M} = (\bff, v)$ where $\bff$ is a pd-frame and $v$ is a valuation for $\bff$. By $\mathcal{M}, w \models \varphi$, we mean that $w \in \tset{\varphi}^\bff(v)$.

Given a pd-frame $\bff = (W, R, B)$, a valuation $v: \Prop \to B$, and a non-empty $U \subseteq W$, the restriction $\bff|_U$ is $(U, R \cap (U \times U), B|_U)$ where $B|_U = \{X \cap U \mid X \in B\}$, and the restriction $v|_U$ is defined by $v|_U(p) = v(p) \cap U$. Then, given $w \in W$ and $n \in \NN$, recall that $R^{\le n}[w]$ is the set of worlds reachable from $w$ in at most $n$ steps (including $w$ itself), and we define the depth-$n$ subframe $\bff_{w, n}$ generated from $w$ by $\bff|_{R^{\le n}[x]}$. Similarly, $v_{w, n} = v|_{R^{\le n}[w]}$. For any pd-model $\mathcal{M} = (\bff, v)$ and $w \in \bff$, $n \in \NN$, $\mathcal{M}_{w, n} = (\bff_{w, n}, v_{w, n})$.


For any $\varphi \in \mathcal{L}$, the modal depth $md(\varphi)$ of $\varphi$ is defined as usual with propositional quantifiers ignored.
\end{definition}

Now we show the standard semantic preservation under passing to point-generated models.
The version without propositional domain has been shown in \cite{ten2006expressivity}.
\begin{lemma}\label{lem:submodel} For any $\varphi \in \mathcal{L}$, $n \ge md(\varphi)$,
    pd-model $\mathcal{M} = (\bff, v)$ and $w \in \bff$, $\mathcal{M}, w \models \varphi$ iff $\mathcal{M}_{w, n}, w \models \varphi$.
\end{lemma}

\begin{proof}
    Fix $\bff = (W, R, B)$, $w \in W$, and $n \in \NN$. Now we prove by induction on formulas that for any $\varphi \in \mathcal{L}$ with $md(\varphi) \le n$, letting $d = md(\varphi)$, for any $x \in \bff_{w, n-d}$ and any valuation $v$ for $\bff$, $x \in \tset{\varphi}^\bff(v)$ iff $x \in \tset{\varphi}^{\bff_{w, n}}(v_{w, n})$. This clearly entails the stated lemma.

    The base case and the Boolean inductive cases are trivial. For the modal case, take a formula $\Diamond\varphi$ such that $md(\Diamond\varphi) \le n$. Let $d = \md(\Diamond\varphi) = \md(\varphi) + 1$. Then $n-md(\varphi) = n-d+1$. Take any $x \in \bff_{w, n-d}$ and valuation $v$ for $\bff$. Then we have the following chain of equivalences where the third equivalence uses the IH:
    \begin{itemize}
        \item $x \in \tset{\Diamond\varphi}^\bff(v)$;
        \item there is $y \in R[x]$ such that $y \in \tset{\varphi}^\bff(v)$;
        \item there is $y \in \bff_{w, n-md(\varphi)} \cap R[x]$ such that $y \in \tset{\varphi}^\bff(v)$;
        \item there is $y \in \bff_{w, n-md(\varphi)} \cap R[x]$ such that $y \in \tset{\varphi}^{\bff_{w, n}}(v_{w, n})$;
        \item there is $y \in R[x]$ such that $y \in \tset{\varphi}^{\bff_{w, n}}(v_{w, n})$;
        \item $x \in \tset{\Diamond\varphi}^{\bff_{w, n}}(v_{w, n})$.
    \end{itemize}

    Now consider the quantifier case. First, suppose $x \in \tset{\exists p\varphi}^\bff(v)$. Then there is $X \in B$ such that $x \in \tset{\varphi}^\bff(v[X/p])$. Apply the IH to $v[X/p]$, and we have that $x \in \tset{\varphi}^{\bff_{w, n}}(v[X/p]_{w, n})$. But note that $v[X/p]_{w, n}$ is just $v_{w, n}[X \cap \bff_{w, n}/p]$, and by the definition of $B_{w, n}$, this is a valuation for $\bff_{w, n}$. So $x \in \tset{\exists p\varphi}^{\bff_{w, n}}(v_{w, n})$. Conversely, suppose $x \in \tset{\exists p\varphi}^{\bff_{w, n}}(v_{w, n})$. Then there is a set $Y \in B_{w, n}$ such that $x \in \tset{\varphi}^{\bff_{w, n}}(v_{w, n}[Y/p])$. Again, by the definition of $B_{w, n}$, there is $X \in B$ such that $Y = X \cap \bff_{w, n}$. Then $v_{w, n}[Y/p]$ is just $v[X/p]_{w, n}$. Apply the IH to $v[X/p]$, then $x \in \tset{\varphi}^{\bff}(v[X/p])$. Then $x \in \tset{\exists p \varphi}^\bff(v)$.
\end{proof}

To separate $p \to \Box p$ from $\mathsf{K}_\Pi\ttt{TMEQAt}^1\ttt{R}^1$, we use the veiled
recession frame $\bff = (\mathbb{Z},R,B)$ that is also used in \cite{van1978two}. Here $R$
is defined by $nRm$ iff $m \ge n-1$ and $B$ is the collection of sets that are `eventually
settled': $X \in B$ iff there is $n \in \mathbb{Z}$ such that for all $m \ge n$, $m \in X$
iff $n \in X$ (we can call this property \emph{settled after $n$}). Using the notation
$[n, \infty)$ for the set of integers $m \ge n$, $B = \{X\subseteq\mathbb{Z} \mid \exists
n \in \mathbb{Z}, [n, \infty) \subseteq X \text{ or } [n, \infty) \cap X = \varnothing\}$.

For a valuation $v:\Prop \to B$ and $p \in \Prop$, we say that $p$ is \emph{settled as
true (by $v$) after $n$}, if $m\in v(p)$ for all $m\geq n$, and $p$ is \emph{settled as
false (by $v$) after $n$}, if $m\not\in v(p)$ for all $m\geq n$.

\begin{lemma}
$(\mathbb{Z},R,B)$ is a quantifiable frame that validates
$\mathtt{T,M,E,Q}$,and $\mathtt{At}^1$ and $\mathtt{R}^1$.
\end{lemma}

\begin{proof}
It is easy to verify that the frame validates \ttt{T}. \ttt{M} is valid since every set in
$B$ is eventually settled. To verify \ttt{E}, notice if a number $n$ satisfies $\Diamond
p\land \Box q$, then all numbers smaller than $n$ would satisfy $\Diamond p$ and numbers
greater would satisfy $\Box q$. To verify \ttt{Q}, notice for arbitrary $nRm$ there is a
finite path from $n$ back to $m$, hence if a world satisfies $\Diamond p\land \Box (p\to
\Box p)$, then it sees a world with $p$ and $p$ can be forced back to itself since all
worlds along the path have $p\to \Box p$. To see that it validates $\mathtt{At}^1$ and
$\mathtt{R}^1$, notice all singleton subsets of $\mathbb{Z}$ are in $B$ and all $R[X]$ are
in $B$ for arbitrary $X\subseteq \mathbb{Z}$.

Now we verify that $\bff = (\mathbb{Z}, R, B)$ is a quantifiable frame. By induction, we
show that for any formula $\varphi \in \mathcal{L}$, for any valuation $v: \Prop \to B$,
$\tset{\varphi}^\bff(v)$ is in $B$. The base case is trivial. As $B$ is clearly closed under
union and complementation, the Boolean cases are also easy. The modal case requires us to
check that $B$ is closed under $m_\Diamond$. For any $X \in B$, if $[n, \infty) \subseteq
X$, then $m_\Diamond(X) = \mathbb{Z}$; and if $[n, \infty) \cap X = \varnothing$, then
$[n+1, \infty) \cap m_\Diamond(X) = \varnothing$. Thus, either way, $m_\Diamond(X) \in B$. 

For the quantifier case, take any existential formula $\exists p \varphi$ and valuation
$v$, and we need to show that $\tset{\exists p \varphi}^\bff(v)$ is in $B$. Let $l$ be an
integer such that for all $q \in \fv(\exists p\varphi)$, $v(q)$ is settled after $l$.
Without loss of generality, we can also assume that for all $q \in \Prop \setminus
\fv(\exists p \varphi)$, $v(q) = \varnothing$. Then for all $q \in \Prop$, $v(q)$ is
settled after $l$. Now we discuss two cases (they are negations of each other):
\begin{itemize}
\item For any $n \ge l+md(\varphi)$ and $X \in B$, $n \not \in \tset{\varphi}^\bff(v[X/p])$.
\item There is $n \ge l + md(\varphi)$ and $X\in B$ such that $n\in \tset{\varphi}^\bff(v[X/p])$.
\end{itemize}
If the former, then $[l + md(\varphi), \infty) \cap \tset{\exists p\varphi}^\bff(v) =
\varnothing$ and thus $\tset{\exists p \varphi}^\bff(v) \in B$. If the latter, then we can
show that $[n, \infty) \subseteq \tset{\exists p \varphi}^\bff(v)$. For any $m \ge n$,
consider the model $\mathcal{M} = (\bff, v[X/p])$ and $\mathcal{M}' = (\bff, v[X+(m-n)/p])$. Note that $B$ is closed under the $+(m-n)$ operation, so $v[X+(m-n)/p]$ is a valuation for $\bff$. Given the case we are discussing, $\mathcal{M}, n \models \varphi$. By Lemma \ref{lem:submodel}, $\mathcal{M}_{n, md(\varphi)}, n \models \varphi$. Now observe that $\mathcal{M}_{n, md(\varphi)}$ and $\mathcal{M}'_{m, md(\varphi)}$ are isomorphic by the isomorphism $\pi: x \mapsto x + (m-n)$. So $\mathcal{M}'_{m, md(\varphi)}, m \models \varphi$. By Lemma \ref{lem:submodel} again, $\mathcal{M}', m \models \varphi$. This means $m \in \tset{\exists p \varphi}^\bff(v)$. As $m$ is chosen arbitrarily from $[n, \infty)$, $[n, \infty) \subseteq \tset{\exists p \varphi}^\bff(v)$.
\end{proof}
Given the above lemma, the set $\Lambda$ of validities of $(\mathbb{Z}, R, B)$ includes $\mathsf{K}_\Pi\ttt{TMEQAt}^1\ttt{R}^1$. But $\Lambda$ does not contain $p \to \Box p$, as it would be false at $0$ if $v(p) = \{0\}$ for example. Thus, $p \to \Box p$ is not in $\mathsf{K}_\Pi\ttt{TMEQAt}^1\ttt{R}^1$, and consequently $\mathsf{K}_\Pi\ttt{TMEQAt}^1\ttt{R}^1$ is not $\ttt{TMEQ}\pi+$.

\end{document}